\documentclass[preprint, 12pt]{elsarticle}
\usepackage[top=1in, bottom=1in, left=1in, right=1in]{geometry}

\usepackage{amsmath,amsfonts,amsthm,amssymb}
\usepackage{graphicx}
\usepackage{verbatim}
\usepackage{color}
\usepackage{amscd}
\usepackage{verbatim}

\newcommand{\loc}{{\mathrm{loc}}}
\newcommand{\core}{C_0^{\infty}(\Omega)}

\newcommand{\be}{\begin{equation}}
\newcommand{\ee}{\end{equation}}
\newcommand{\mysection}[1]{\section{#1}\setcounter{equation}{0}}
\newcommand{\bea}{\begin{eqnarray}}
\newcommand{\eea}{\end{eqnarray}}
\newcommand{\bean}{\begin{eqnarray*}}
\newcommand{\eean}{\end{eqnarray*}}

\newcommand{\R}{{\mathbb R}}

\newtheorem{theorem}{Theorem}[section]
\newtheorem{corollary}[theorem]{Corollary}
\newtheorem{lemma}[theorem]{Lemma}

\newtheorem{definition}[theorem]{Definition}
\newtheorem{remark}[theorem]{Remark}
\newtheorem{proposition}[theorem]{Proposition}

\newtheorem{problem}[theorem]{Problem}
\newtheorem{problems}[theorem]{Problems}
\newtheorem{example}[theorem]{Example}

\newtheorem{Rem}[theorem]{Remark}

\newtheorem{Cor}[theorem]{Corollary}

\newtheorem*{corollary*}{Corollary}
\newtheorem*{theorem*}{Theorem}
\newtheorem*{definition*}{Definition}

\catcode`@=11 \@addtoreset{equation}{section} \catcode`@=12
\newcommand{\Real}{\mathbb{R}}

\newcommand{\Nat}{\mathbb{N}}

\newlength{\wex}  \newlength{\hex}
\newcommand{\dx}{\,\mathrm{d}x}
\newcommand{\dy}{\,\mathrm{d}y}

\def\ga{\alpha}     \def\gb{\beta}       \def\gg{\gamma}
       \def\gd{\delta}      
\def\gth{\theta}                         \def\vge{\varepsilon}
\def\gf{\phi}       \def\vgf{\varphi}    
            \def\gl{\lambda}

\def\gs{\sigma}       
      \def\gw{\omega}
                
     \def\Gd{\Delta}      \def\Gf{\Phi}

\def\Gw{\Omega}              
\begin{document}

\begin{frontmatter}
\title{Criticality theory of half-linear equations with the $(p,A)$-Laplacian}

\author{Yehuda Pinchover}
\address{Department of Mathematics, Technion - Israel Institute of
Technology,  Haifa, Israel}
 \ead{pincho@techunix.technion.ac.il}
\author{Netanel Regev}
\address{Department of Mathematics, Technion - Israel Institute of
Technology,  Haifa, Israel}
\ead{nati@techunix.technion.ac.il}

\begin{abstract}
We study positive solutions of half-linear second-order elliptic equations of the form
$$Q_{A,V}(u):=  -\mathrm{div} (|\nabla u|_{A}^{p-2}A(x)\nabla u)+ V(x)|u|^{p-2}u=0 \quad \mbox{in }\Gw,$$
where $1<p<\infty$, $\Omega$ is a domain in $\mathbb{R}^{n}$, $n\geq 2$, $V\in L_{\mathrm{loc}}^{\infty}(\Omega)$, $A=\big(a_{ij}\big)\in L_{\mathrm{loc}}^{\infty}(\Omega,\R^{n^2})$ is a symmetric and locally uniformly positive definite matrix  in $\Omega$,  and $$|\xi|_{A}^{2}:=\left\langle A(x)\xi,\xi\right\rangle=\sum_{i,j=1}^n a_{ij}(x)\xi_i\xi_j \qquad x\in \Gw, \xi=(\xi_1,\ldots,\xi_n)\in\R^n.$$
We extend criticality theory which has been established for linear operators and for half-linear operators involving the $p$-Laplacian, to the operator $Q_{A,V}$. We prove Liouville-type theorems, and study the behavior of positive solutions of the equation $Q_{A,V}(u)=0$ near an isolated singularity and near infinity in $\Gw$, and obtain some perturbations results.

\medskip

\noindent {\em Keywords.}  Half-linear equations, Liouville-comparison principle, maximum principles, $p$-Laplacian, positive solutions, quasilinear elliptic
equation, removable singularity.
%%%%%%%%%%%%%%%%
%\date{}
%\thanks{}
%\subjclass[]{}
%\keywords{}
%%%%%%%%%%%%%%%
\end{abstract}
\end{frontmatter}

%\date{}
%\thanks{}
%\subjclass[]{}
%\keywords{}

%\maketitle
\mysection{Introduction}\label{sec1}
In this paper we study positive solutions of half-linear elliptic partial differential equations of second-order.
Recall that a partial differential equation $Q(u)=0$ is said to be {\em half-linear} if for any $\ga\in \R$ we have $Q(\ga v)=0$, whenever $Q(v)=0$. So, half-linear equations satisfy the homogeneity property of linear equations but not the additivity. Therefore, it is natural to expect that positive solutions of such equations would share some fundamental properties of positive solutions of linear elliptic equations \cite[and references therein]{P07}. It turns out that this is indeed the case for certain half-linear equations. In fact, the theory of positive solutions of
half-linear elliptic equations associated with the $p$-Laplacian operator $\Gd_p$ and a potential term $V$ has been
studied extensively in recent years (see for example,
\cite{{AH1},{AH2},{GS},{PTT},{PT1},{PT2},{PT3}} and the references therein). In particular, the criticality theory and especially the well
known Agmon-Allegretto-Piepenbrink (AAP) theorem has been extended from the linear case to such half-linear equations (see \cite[Theorem 2.12]{CFKS} and \cite[Theorem~2.3]{PT1}).

In the present work we extend some positivity results of the above mentioned papers concerning half-linear equations with the $p$-Laplace operator, to the case of half-linear equations with the so called $(p,A)$-Laplace operator, where $A$ is a given matrix. More precisely, the above mentioned papers study positivity properties of the functional
\[
\mathcal{Q}_{V}(\vgf):=\int_{\Omega}(|\nabla \vgf|^{p}+V(x)|\vgf|^{p})\dx \qquad \vgf\in \core,
\]
and its associated Euler-Lagrange equation
\begin{equation}\label{derivative0}
Q_V(u):= \frac{1}{p}\mathcal{Q}_V'(u)=-\Gd_p(u)+V(x)|u|^{p-2}u=0 \quad \mbox{in } \Gw,
\end{equation}
where $1\!<p<\!\infty$, $\Omega$ is a domain in $\mathbb{R}^{n}$, $n\!\geq\! 2$. $\Gd_p(u)\!:=\!\nabla\cdot(|\nabla u|^{p-2}\nabla u)$ is the celebrated {\em $p$-Laplacian}, and $V\in L_{\mathrm{loc}}^{\infty}(\Omega)$ is a real potential.

In our work we study the functional
\[
\mathcal{Q}_{A,V}(\vgf):=\int_{\Omega}(|\nabla \vgf|_{A}^{p}+V(x)|\vgf|^{p})\dx\qquad \vgf\in C_{0}^{\infty}(\Omega),
\]
and its Euler-Lagrange equation
\begin{equation}\label{derivative}
Q_{A,V}(u):=\frac{1}{p}\mathcal{Q}'_{A,V}(u) =  -\nabla \cdot (|\nabla u|_{A}^{p-2}A(x)\nabla u)+ V(x)|u|^{p-2}u=0 \quad \mbox{in } \Gw,
\end{equation}
where $p$, $\Gw$, and $V$ are as assumed above, $A=\big(a_{ij}\big)\in L_{\mathrm{loc}}^{\infty}(\Omega, \R^{n^2})$ is a symmetric and locally uniformly positive definite matrix  in $\Omega$, and $$|\xi|_{A}^{2}:=\left\langle A(x)\xi,\xi\right\rangle=\sum_{i,j=1}^n a_{ij}(x)\xi_i\xi_j \qquad x\in \Gw, \xi=(\xi_1,\ldots,\xi_n)\in\R^n.$$

The aim of the paper is to establish criticality theory for the operator $Q_{A,V}$. In particular, we prove Liouville-type theorems, study the behavior of positive solutions of equation \eqref{derivative} near an isolated singularity and near infinity in $\Gw$, and obtain some perturbations results.

It is worth noting that in \cite{HKM} the authors study, under stronger assumptions on the matrix $A$, the case where the potential term is missing (i.e. $V=0$). We generalize some results of this monograph to the case where $V \neq 0$ without assuming \cite[(3.4) and (3.5)]{HKM}. On the other hand, throughout the paper we always assume that $V\in L^\infty_\loc(\Gw)$. It would be interesting to extend our results to the case where $V\in L^q_\loc(\Gw)$ for an appropriate $q\geq 1$.

The outline of the present paper is as follows. Section~\ref{sec_Pre} is devoted to some preliminaries, while in Section~\ref{sec_picone} we extend to our setting some fundamental tools, like the {\em Picone identity} \cite{{AH1},{AH2},Y}, the {\em Anane-D\'iaz-Saa identity} \cite{DS,AA}, and the {\em simplified energy} \cite{PTT}.

Section~\ref{sec_SMP} is mainly devoted to generalizations of results of J.~Garc\'ia-Meli\'an and J.~Sabina de Lis \cite{GS} concerning the relationships between the principal eigenvalue, the weak and strong maximum principles, and the solvability of the Dirichlet problem.
We note that under the assumptions of \cite{GS}, solutions of \eqref{derivative0} are $C^{1,\ga}$-smooth and satisfy the boundary point lemma. On the other hand, under our assumptions, solutions of \eqref{derivative} are only $C^{\ga}$, and therefore, we need to provide completely new proofs to some of the results in \cite{GS}.

In Section~\ref{sec_AAP} we extend to our case the well-known AAP theorem dealing with the relationships between positivity properties of the functional $\mathcal{Q}_{A,V}$ and the existence of positive (super)solution of the equation $Q_{A,V}(u)=0$.
Section~\ref{sec_main_thm} contains the proof of the Main Theorem (Theorem~\ref{main thm}) which generalizes \cite[Theorem 3.3]{PT2} and concerns characterizations of critical/subcritical operators. Using the Main Theorem, we study in Section~\ref{sec_crit_theory} criticality properties of the functional $\mathcal{Q}_{A,V}$, generalizing results of \cite{{PT1},{PT2},{PT3}}.

In Section~\ref{sec_Liouville}  we use the simplified energy to generalize the Liouville-comparison principle proved in \cite[Theorem~1.9]{PTT}, while
in Section~\ref{sec_minimal growth} we prove the existence of a positive minimal Green functions $G_{A,V}^\Gw (x,x_0)$ in the subcritical case, and global minimal solutions of the equation $Q_{A,V}(u)=0$ in $\Omega$ in the critical case. The problem concerning the uniqueness of $G_{A,V}^\Gw (x,x_0)$ remains open for the case $A\neq I$.
%%%%%%%%%%%%%%%%%%%%%%%%%%%%%%%%%%%%%%%%%%%%%%%%%%%
\section{Preliminaries} \label{sec_Pre}
%%%%%%%%%%%%%%%%%%%%%%%%%%%%%%%%%%%%%%%%%%%%%%%%%%%%%
In this section we fix our setting and notations, and introduce some basic definitions. Throughout the paper $1<p<\infty$, and $\Omega \subseteq \Real^{n}$ is a domain.
We write $S \Subset \Omega$ if $\Omega$ is open, $\overline{S}$ is
compact and $\overline{S} \subset \Omega$. By an {\em exhaustion} of $\Omega$ we mean a sequence $\left\{\Omega_{N}\right\}$ of smooth, relatively compact domains such that $x_{0} \in \Omega_{1}$, $\Omega_{N} \Subset \Omega_{N+1}$, and $\bigcup_{N=1}^{\infty}\Omega_{N}=\Omega$.

Let $f,g:\Gw \to [0,\infty)$. We denote $f \asymp g$ if there exists a positive constant $C$ such that $C^{-1}g \leq f \leq Cg$ in $\Gw$. Also, $f\gneqq g$ if $f\geq 0$ but $f\neq 0$. For $1< p< \infty$ we denote $p\ \!':=p/(p-1)$ the conjugate index of $p$.
$B_r(x)$ is the open ball of radius $r$ centered at $x$.

\medskip

Let us present the regularity assumptions for the operator $Q_{A,V}$ which ensure the validity of the weak and strong Harnack inequalities, and the $C^\ga$-regularity of solutions.  Throughout our paper we assume (unless otherwise stated) that
\begin{align*}\label{assumption AV}
\tag{A} & A:\Gw \to \R^{n^2} \mbox{ is a symmetric measurable matrix}.  \\[4mm]
\tag{E} &  \forall K\Subset \Gw \; \exists\, \theta_K>0\; \mbox{ s.t. }  \theta_K|\xi|^{2} \leq \left\langle A(x)\xi,\xi\right\rangle
\leq\theta_K^{-1} |\xi|^2\;\;\; \forall x\in K, \xi \in \Real^{n}.\\[4mm]
\tag{V} & V \in L_{\mathrm{loc}}^{\infty}(\Omega).
\end{align*}
\begin{Rem}{\em
Some results of the paper are proved under stronger regularity assumptions. Indeed, for the up to the boundary $C^{1,\ga}$-regularity we need to assume that $A\in C^{\ga}$, while for the validity of the boundary point lemma (Corollary~\ref{corSMP}) we need to assume
that $A\in C^{2}$.
}
\end{Rem}
We now introduce a formal differential operator $\Delta_{p,A}(u):=\nabla \cdot (|\nabla u|_{A}^{p-2}A(x)\nabla u)$,  called the {\em $(p,A)$-Laplacian}.
For $u \in W_{\mathrm{loc}}^{1,p}(\Omega)$  we define:
\begin{equation}\label{}
\left\langle -\Delta_{p,A}(u),\vgf\right\rangle_\Gw:=\int_{\Omega}|\nabla u|_{A}^{p-2}A(x)\nabla u\cdot \nabla \vgf\dx \qquad \forall \vgf \in \core.
\end{equation}\label{eq_LapAV}
\begin{remark}\label{rem_LapAV}{\em
Clearly, for $\Gw'\Subset \Gw$ and $u,v\in W^{1,p}(\Omega')$, the integral $\int_{\Omega'}|\nabla u|_{A}^{p-2}A(x)\nabla u \cdot\nabla v\dx $
is well defined, and for such functions we shall still denote
$$\langle -\Delta_{p,A}(u),v \rangle_{\Gw'}:=\int_{\Omega'}|\nabla u|_{A}^{p-2}A(x)\nabla u\cdot \nabla v\dx.$$
 }
\end{remark}

\begin{definition}
{\em A function $v \in W_{\loc}^{1,p}(\Omega)$ is a {\em(weak) solution} of equation \eqref{derivative}
if
\begin{equation}
\label{sol}
\left\langle Q_{A,V}(v),\varphi\right\rangle_{\Omega}:=\int_{\Omega}{(|\nabla v|_{A}^{p-2}A(x)\nabla v \cdot \nabla \varphi+V(x)|v|^{p-2}v\varphi)\dx}=0 \qquad \forall\varphi \in C_{0}^{\infty}(\Omega).
\end{equation}

We say that a positive function $v \in  W_{\loc}^{1,p}(\Omega)$ is a {\em supersolution (resp., subsolution)} of equation \eqref{derivative}
if for every nonnegative $\varphi \in C_{0}^{\infty}(\Omega)$
\begin{equation}
\label{super}
\int_{\Omega}{(|\nabla v|_{A}^{p-2}A(x)\nabla v \cdot \nabla \varphi+V|v|^{p-2}v\varphi)\dx}\geq0.\ \ (\mbox{resp.} \leq0).
\end{equation}

By {\em uniqueness of positive (super)solutions} of \eqref{derivative}, we always mean uniqueness up to a multiplicative constant.
 }
\end{definition}

%%%%%%%%%%%%%%%%%%%%%%%%%%%%%%%%%%
In the sequel we need the following elementary lemma.
\begin{lemma}\label{subsol}
Let $v\in W^{1,p}_{\mathrm{loc}}(\Omega)$ be a subsolution  of the
equation \eqref{derivative}, and let $v_+(x):=\max\{0,v(x)\}$. Then $v_+$ is also a subsolution of \eqref{derivative}.
\end{lemma}
\begin{proof} The proof is a slight modification of the proof of \cite[Lemma~2.4]{PTT} obtained  by replacing the Euclidean inner product with the inner product $\left\langle A(x)\xi,\eta\right\rangle$ induced by the matrix $A$.
\end{proof}
%%%%%%%%%%%%%%%%%%

\begin{definition}\label{def_pos} {\em
The functional $\mathcal{Q}_{A,V}$ is {\em nonnegative} in $\Omega$ (notation: $\mathcal{Q}_{A,V}\geq 0$  in $\Omega$) if
\begin{equation} \label{func}
\mathcal{Q}_{A,V}(\vgf)=\int_{\Omega}(|\nabla \vgf|_{A}^{p}+V(x)|\vgf|^{p})\dx \geq 0 \qquad \vgf\in C_{0}^{\infty}(\Omega).
\end{equation}
}
\end{definition}
% Throughout the present paper we assume that $\mathcal{Q}_{A,V}$ is {\em nonnegative} in $\Omega$.

\begin{definition}{\em

Assume that $\mathcal{Q}_{A,V}\geq 0$ in $\Gw$. We say that $\mathcal{Q}_{A,V}$ is {\em subcritical} in $\Omega$ if there exists a nonzero nonnegative continuous function $W$ in $\Omega$ such that
\begin{equation}\label{subcr}
\mathcal{Q}_{A,V}(\vgf) \geq \int_{\Omega}{W|\vgf|^{p}\dx}  \qquad \forall \vgf \in C_{0}^{\infty}(\Omega).
\end{equation}
A nonnegative functional $\mathcal{Q}_{A,V}$  in $\Omega$ which is not subcritical in $\Gw$ is called {\em critical} in $\Gw$.
A functional $\mathcal{Q}_{A,V}$ is {\em supercritical} in $\Omega$ if $\mathcal{Q}_{A,V}$ is not nonnegative in $\Gw$.

}
\end{definition}

\begin{definition}{\em
A sequence $\{\vgf_{k}\}\!\subset\! C_{0}^{\infty}(\Omega)$ is a {\em null sequence} with respect to the nonnegative functional $\mathcal{Q}_{A,V}$ in $\Gw$ if $\vgf_{k} \!\geq \!0$ for all $k \in \Nat$, and there exists an open set $B \Subset \Omega$  such that
\begin{equation}
\lim_{k \to \infty}\mathcal{Q}_{A,V}(\vgf_k)=\lim_{k \to \infty}\int_{\Omega}(|\nabla \vgf_{k}|_{A}^{p}+V|\vgf_{k}|^{p})\dx=0, \quad \mbox{and } \int_{B}{|\vgf_{k}|^{p}\dx}=1.
\end{equation}
We say that a positive function $\phi\in W_{\loc}^{1,p}(\Omega)$ is  {\em Agmon's ground state} (or simply a ground state) of the functional $\mathcal{Q}_{A,V}$ in $\Omega$ if $\phi$ is an $L_{\loc}^{p}(\Omega)$ limit of a null sequence.
}
\end{definition}
%%%%%%%%%%%%%%%%%%%%%%
\begin{Rem}\label{ground} {\em
The requirement that $\left\{\varphi_{k}\right\} \subset C^{\infty}_{0}(\Omega)$, can be weakened by assuming only that $\left\{\varphi_{k}\right\} \subset W^{1,p}_{0}(\Omega)$, and that $\int_{B}|\varphi_{k}|^{p}\dx \asymp 1$ (instead of $\int_{B}|\varphi_{k}|^{p}\dx=1$).
}
\end{Rem}

\begin{example}[\cite{PTT} Example~1.7]\label{ex_subcr}{\em
Let $A\in C(\R^N,\R^{n^2})$ be a symmetric, bounded, and uniformly positive definite matrix in $\R^n$. Consider the functional
$$\mathcal{Q}_{A,0}(\varphi):=\int_{\Real^{n}}|\nabla \varphi|_A^{p}\dx\qquad \varphi\in C^\infty_0(\R^n).$$
If $p\geq n$, then \cite[Theorem~2]{MP} and Theorem~\ref{main thm} imply that $\mathcal{Q}_{A,0}$ is critical in $\R^n$ and that $\phi=\mathrm{const.}>0$ is its  ground state (see Example~\ref{ex2} for an extension of this result).

On the other hand, if $p<n$, then the equation $Q_{I,0}(u)=-\Delta_{p}(u)=0$ in $\R^n$
admits two linearly independent positive supersolutions
 $$u(x):=\mathrm{const.}, \qquad \mathrm{and} \qquad v(x):=\left[1+|x|^{\frac{p}{p-1}}\right]^{\frac{p-n}{p}}.$$ Hence, Theorem~\ref{main thm} implies that $-\Delta_{p}$ is subcritical in $\Real^{n}$. For further examples see \cite{PTT}.
 }
\end{example}
We conclude the present section with the following well known compactness result.

\noindent\textbf{Harnack convergence principle}. Let $\left\{\Omega_{N}\right\}$ be an exhaustion of $\Omega$. Assume that $\left\{A_{N}\right\}^{\infty}_{N=1}$ is a sequence of symmetric and  positive definite matrices satisfying $A_N\in L^{\infty}(\Omega_N,\R^{n^2})$ such that the sequence $\left\{A_{N}\right\}^{\infty}_{N=1}$ converges locally uniformly to a matrix $A$ satisfying conditions (A) and (E). Assume also that  $V_{N} \in L^{\infty}(\Omega_N)$ satisfy $V_{N} \to V$ in $L^{\infty}_{\loc}(\Omega)$. For each $N\geq 1$, let $u_{N}$ be a positive solution of the equation
\begin{equation}
Q_{A_{N},V_{N}}(w)=:-\Delta_{A_{N},p}(w)+V_{N}|w|^{p-2}w=0\qquad \mathrm{in} \ \Omega_{N},
\end{equation}
satisfying $u_{N}(x_{0})=1$.
By Harnack's inequality and elliptic regularity \cite[Chapter~7]{PS}, and a diagonalization argument, there exist $0< \beta <1$, and a subsequence $\left\{u_{N_{k}}\right\}$ of $\left\{u_{N}\right\}$ that converges in $C^{\beta}_{\loc}(\Omega)$ to a positive solution $u\in W^{1,p}_{\mathrm{loc}}(\Gw)$ of the equation  $$Q_{A,V}(w)=-\Delta_{A,p}(w)+V|w|^{p-2}w=0 \qquad \mbox{in } \Omega.$$

%%%%%%%%%%%%%%%%%%%%%%%%%%%%%%%%%%%%%%%%%%%%%%%%%%%
\mysection{Picone identity}\label{sec_picone}
%%%%%%%%%%%%%%%%%%%%%%%%%%%%%%%%%%%%%%%%%%%%%%%%%%%%%
We start with a simple generalization of Picone identity (cf.~\cite{{AH1},{AH2},{Y}}).

\begin{proposition}
\label{Picone}
Let $v>0$, $u\geq0$ be differentiable functions in $\Gw$, and let $A$ satisfy assumptions (A) and (E) in $\Gw$.
Denote
\begin{equation}\label{L_A}
L_A(u,v)(x):=|\nabla u|_{A}^{p}+(p-1)\frac{u^{p}}{v^{p}}|\nabla v|_{A}^{p}-p\frac{u^{p-1}}{v^{p-1}}\nabla u\cdot A(x)\nabla v|\nabla v|_{A}^{p-2},
\end{equation}
and
\begin{equation}
R_A(u,v)(x):=|\nabla u|_{A}^{p}-\nabla \Big(\frac{u^{p}}{v^{p-1}}\Big)\cdot A(x)\nabla v|\nabla v|_{A}^{p-2}.
\end{equation}
Then
\begin{equation} \label{L=R}
L_A(u,v)(x)=R_A(u,v)(x)\geq 0 \qquad \forall x\in \Gw.
\end{equation}
Moreover, $L_A(u,v)=0$ a.e. in $\Omega$ if and only if $u=kv$  in $\Omega$ for some constant $k\geq 0$.

\end{proposition}

\begin{proof} Use the proof of \cite[Theorem~1.1]{AH1} and just replace the Euclidean inner product $\left\langle \xi,\eta\right\rangle$ with the inner product $\left\langle A(x)\xi,\eta\right\rangle$ induced by the matrix $A$.
\end{proof}
\begin{remark} \label{rem_PI}{\em
By a standard approximation argument, it follows that Proposition~\ref{Picone} holds true if $v>0$, $u\geq0$ are in $W^{1,p}_{\mathrm{loc}}(\Omega)$ and
$uv^{-1}\in L^\infty_{\mathrm{loc}}(\Omega)$.
 }
\end{remark}

\begin{proposition} \label{Q=L}
Let $A$ and $V$ satisfy assumptions (A), (E) and (V).  Let $v\in W_{\loc}^{1,p}(\Omega)$, be a positive solution (\mbox{resp.} supersolution) of (\ref{derivative}). Then for every $u\geq 0$, $u\in W_{\loc}^{1,p}(\Omega)$ with compact support in $\Gw$  such that $\frac{u}{v}\in L^\infty_{\mathrm{loc}}(\Gw)$ we have

\begin{equation}\label{eq_Q=L}
    \mathcal{Q}_{A,V}(u)\!=\!\int_{\Omega}\!L_{A}(u,v)\dx\geq 0. \qquad\left(\mbox{resp. } \; \mathcal{Q}_{A,V}(u)\!\geq\! \int_{\Omega}\!L_{A}(u,v)\dx \!\geq\! 0 \right).
\end{equation}
Moreover, \eqref{eq_Q=L} holds true if $\Omega\Subset \Real^{n}$, $A$ is a bounded measurable symmetric matrix which is uniformly positive definite in $\Omega$, $V\in L^\infty(\Gw)$,
$v\in W^{1,p}(\Omega)$ is a positive solution (\mbox{resp.} supersolution) of (\ref{derivative}), and $0\leq u\in W_0^{1,p}(\Omega)$ such that $\frac{u}{v}\in L^\infty(\Gw)$.
\end{proposition}

\begin{proof}
From (\ref{L=R}) we have that
\begin{equation}\label{eq_QAV_geq_LA}
0\leq \int_{\Omega}L_A(u,v)\dx=\int_{\Omega}{|\nabla u|_{A}^{p}}\dx-\int_{\Omega}{\nabla \left(\frac{u^{p}}{v^{p-1}}\right)\cdot A(x)\nabla v|\nabla v|_{A}^{p-2}}\dx.
\end{equation}
Using  (\ref{sol}) (resp. (\ref{super})) with $ \frac{u^{p}}{v^{p-1}} \in W_{0}^{1,p}(\Omega)$, we get that
$$ \int_{\Omega}{|\nabla v|_{A}^{p-2}A(x)\nabla v \cdot \nabla \left(\frac{u^{p}}{v^{p-1}}\right)+V|v|^{p-2}v\frac{u^{p}}{v^{p-1}}\dx}=0.\  \left(\mbox{resp. } \geq 0\right), $$
which imply
\begin{equation}\label{eqref}
    -\int_{\Omega}{\nabla \left(\frac{u^{p}}{v^{p-1}}\right)\cdot A(x)\nabla v|\nabla v|_{A}^{p-2}}\dx=\int_{\Omega}V|u|^{p}\dx.\  \left(\mbox{resp. } \leq \int_{\Omega}V|u|^{p}\dx\right).
\end{equation}

Combining \eqref{eq_QAV_geq_LA} with \eqref{eqref}, we arrive at the desired conclusion.
\end{proof}

Let the hypothesis of Proposition~\ref{Q=L} be satisfied. Denote  $w:=\frac{u}{v}$. Then by \eqref{eq_Q=L},
\begin{equation} \label{1}
\mathcal{Q}_{A,V}(vw)\displaystyle{\substack{= \\[1mm] (\geq)}} \int_{\Omega}\Big[|v\nabla w+w \nabla v|_{A}^{p}-w^{p}|\nabla v|^{p}_{A}-pw^{p-1}v|\nabla v|^{p-2}_{A}\nabla v\cdot A(x) \nabla w \Big]\dx .
\end{equation}
Similarly, for a such nonnegative {\em subsolution} $v$ of (\ref{derivative}) we have,
\begin{equation} \label{2}
\mathcal{Q}_{A,V}(vw) \leq
 \int_{\Omega \cap \{v>0\}}\!\!\Big[|v\nabla w+w \nabla v|_{A}^{p}-w^{p}|\nabla v|^{p}_{A}-pw^{p-1}v|\nabla v|^{p-2}_{A}\nabla v\cdot A(x) \nabla w \Big]\dx.
\end{equation}

\medskip

In Proposition~\ref{Q=L}, the functionals $\mathcal{Q}_{A,V}$ and $\int L_{A}\dx$ are both nonnegative, but in general, both functionals contain indefinite terms. Therefore, we show now, as in \cite[Lemma 2.2]{PTT}, that $\mathcal{Q}_{A,V}$ is equivalent to a {\itshape simplified energy} containing only nonnegative terms.

\begin{lemma}\label{lem_simplified}
Let the hypothesis of Proposition~\ref{Q=L} be satisfied, where $v$ is a positive solution of (\ref{derivative}). Let $w:=\frac{u}{v}$. Then
\begin{equation} \label{3}
\mathcal{Q}_{A,V}(vw) \asymp \int_{\Omega}v^{2}|\nabla w|^{2}_{A}\Big(w|\nabla v|_{A}+v|\nabla w|_{A}\Big)^{p-2}\dx.
\end{equation}
If $v \in W^{1,p}_{\loc}(\Omega)$ is a such nonnegative subsolution of (\ref{derivative}), then
\begin{equation} \label{5}
\mathcal{Q}_{A,V}(vw) \leq C\int_{\Omega \cap \{v>0\}}v^{2}|\nabla w|^{2}_{A}\Big(w|\nabla v|_{A}+v|\nabla w|_{A}\Big)^{p-2}\dx.
\end{equation}
\end{lemma}

\begin{proof} We claim that the following estimate holds true for all $x\in \Gw$
\begin{equation} \label{6}
|a+b|^{p}_{A}-|a|^{p}_{A}-p|a|^{p-2}_{A} a\cdot A(x)b \asymp |b|^{2}_{A}\left(|a|_{A}+|b|_{A}\right)^{p-2} \qquad \forall a,b \in \Real^{n},
\end{equation}
where the equivalence constant does not depend on $A(x)$. The proof of \eqref{6} is similar to the proof of \cite[Inequality (2.19)]{PTT}.
Set now $a:=w|\nabla v|_{A}, b:=v|\nabla w|_{A}$, and we obtain (\ref{3}) and (\ref{5}) by applying (\ref{6}) to (\ref{1}) and (\ref{2}), respectively.
\end{proof}

The following lemma is a simple generalization of \cite[Lemma~4]{GS}. It was proved by D\'iaz and Saa in \cite[Lemma~2]{DS} and by Anane in \cite[Proposition~1]{AA}, for the case $A=I$ and $V=0$.

\begin{lemma}
\label{I(u,v)}
Let $\Omega\subset \Real^{n}$ be a bounded domain, and assume that $A$ is a symmetric matrix of class $L^\infty(\Gw)$ which is positive definite in $\Omega$, and $V\in L^\infty(\Gw)$.
Consider the functional
\begin{equation*}
    I_{A}(u,v):=\left\langle Q_{A,V}(u),\frac{u^p-v^p}{u^{p-1}}\right\rangle_{\Omega}-\left\langle Q_{A,V}(v),\frac{u^p-v^p}{v^{p-1}}\right\rangle_{\Omega},
\end{equation*}
defined on $$\mathcal{D}(I_A)=\left\{(u,v)\in  (W^{1,p}(\Omega))^{2} \mid u,v\geq 0 \ \ in\  \ \Omega,\  \frac{u}{v}, \frac{v}{u} \in L^{\infty}(\Omega)\right\},$$
(see Remark~\ref{rem_LapAV}).
Then $I_{A}\geq 0$ on $\mathcal{D}(I_A)$, and $I_{A}(u,v)=0$ if and only if $u=kv$ for some constant $k> 0$.
\begin{proof}
First we note that
\begin{equation*}
\left\langle -\Delta_{p,A}(u),\frac{u^p-v^p}{u^{p-1}}\right\rangle_{\Omega} =
     \int_{\Omega}|\nabla u|_{A}^{p}\dx- \int_{\Omega}|\nabla u|_{A}^{p-2}A(x)\nabla u \cdot\nabla \Big(\frac{v^p}{u^{p-1}}\Big)\dx.                                                    \nonumber
\end{equation*}
So, we have
\begin{multline*}
I_{A}(u,v) =\left\langle Q_{A,V}(u),\frac{u^p-v^p}{u^{p-1}}\right\rangle_{\Omega}-\left\langle Q_{A,V}(v),\frac{u^p-v^p}{v^{p-1}}\right\rangle_{\Omega}=\\
    \left\langle -\Delta_{p,A}(u),\frac{u^p-v^p}{u^{p-1}}\right\rangle_{\Omega}-\left\langle -\Delta_{p,A}(v),\frac{u^p-v^p}{v^{p-1}}\right\rangle_{\Omega}=\\
    \int_{\Omega}\!\!\left(|\nabla u|_{A}^{p}\!-\!|\nabla u|_{A}^{p\!-\!2}A(x)\nabla u \!\cdot\!\nabla \!\Big(\frac{v^p}{u^{p\!-\!1}}\Big)\right)\dx \! +\!
    \int_{\Omega}\!\!\left(|\nabla v|_{A}^{p}\!-\!|\nabla v|_{A}^{p-2}A(x)\nabla v \!\cdot\!\nabla \!\Big(\frac{u^p}{v^{p\!-\!1}}\Big)\right)\dx\!=\\
    \int_{\Omega}R_{A}(u,v)\dx +\int_{\Omega}R_{A}(v,u)\dx,         \nonumber
\end{multline*}
and by Proposition~\ref{Picone} (see Remark~\ref{rem_PI}) we obtain the required result.
\end{proof}

\end{lemma}

\mysection {Maximum principles, the principal eigenvalue and the comparison principle}\label{sec_SMP}
Throughout the present section we assume  that $\Omega \subset \Real ^{n}$ is a {\em bounded} domain, and that the coefficients of the operator $Q_{A,V}$ are bounded,  and $A(x)=\big(a_{ij}(x)\big)$ is a symmetric matrix which is positive definite in $\Omega$ such that for some $0<\theta \leq \Theta$ we have
\begin{equation}\label{eq_theta}
\theta|\xi|^{2} \leq \sum_{i,j=1}^{n}a_{ij}(x)\xi_{i}\xi_{j} \leq\Theta |\xi|^2\qquad \forall x\in \Omega, \xi \in \Real^{n}.
\end{equation}

Let $f \in W^{-1,p'}(\Omega)$.  A function $v \in W_{\loc}^{1,p}(\Gw)$ is a weak solution of the equation $Q_{A,V}(u)=f$ in $\Gw$ if
\begin {equation} \label {sol non}
\int_{\Omega}{(|\nabla v|_{A}^{p-2}A(x)\nabla v \cdot \nabla \varphi+V|v|^{p-2}v\varphi)\dx}=\left\langle f,\varphi\right\rangle \qquad \forall \varphi \in C_{0}^{\infty}(\Omega).
\end {equation}
 For $f\geq 0$, a supersolution of the equation $Q_{A,V}(u)=f$ is defined in a similar fashion.

 Denote by
$$a^{ij}(A,x,\nabla u):= \frac{\partial}{\partial u_{x_j}}\left(|\nabla u|^{p-2}_{A}A_i(x)\cdot \nabla u\right),$$
where  $A_{i}(x)$ denotes the $i$th-row of the matrix $A(x)=\big(a_{ij}(x)\big)$.

\begin{proposition} \label{bound}
Let $A$ satisfy (A) and \eqref{eq_theta}. Then for all $\xi \in \Real^{n}$ and $x\in \Gw$
\begin{equation}\label {positive definite}
\min\{1,p-1\}\theta^{\frac{p}{2}} |\nabla u|^{p-2}|\xi|^{2} \leq \sum_{i,j=1}^{n}a^{ij}(A,x,\nabla u)\xi_{i}\xi_{j}\leq \max\{1,p-1\}\Theta^{\frac{p}{2}} |\nabla u|^{p-2}|\xi|^{2}.
\end{equation}
\end{proposition}

\begin{proof}
A direct calculation shows that
\[
\Big(a^{ij}(A,x,\nabla u)\Big)=|\nabla u|^{p-2}_{A}\Big(a_{ij}(x)+\frac{p-2}{|\nabla u|^{2}_{A}}\left(A_{i}(x)\nabla u\right)\otimes\left(A_{j}(x)\nabla u\right)\Big).
\]

Consequently, for all $\xi \in \Real^{n}$ we obtain that
\begin{equation}\label{eq4_4}
\sum_{i,j=1}^{n}a^{ij}(A,x,\nabla u)\xi_{i}\xi_{j}=  |\nabla u|^{p-2}_{A}\Big[\xi^{T}A(x)\xi+\frac{p-2}{|\nabla u|^{2}_{A}}\left\langle A(x)\nabla u,\xi\right\rangle^{2}\Big].
\end{equation}
Obviously, for $p \geq 2$ we have that $$\sum_{i,j=1}^{n}a^{ij}(A,x,\nabla u)\xi_{i}\xi_{j}\geq \theta|\nabla u|^{p-2}_{A}|\xi|^{2}.$$
For $1<p<2$, \eqref{eq4_4} and the Cauchy-Schwarz inequality imply
\begin{equation*}
\sum_{i,j=1}^{n}a^{ij}(A,x,\nabla u)\xi_{i}\xi_{j} \geq
  |\nabla u|^{p-2}_{A}\Big[|\xi|^{2}_{A}+\frac{p-2}{|\nabla u|^{2}_{A}}|\nabla u|^{2}_{A}|\xi|^{2}_{A}\Big] =(p-1)\,\gth|\nabla u|^{p-2}_{A}|\xi|^{2}.
\end{equation*}

For the upper bound, we see that (\ref{eq4_4}) readily implies for $1<p\leq2$ that
 $$\sum_{i,j=1}^{n}a^{ij}(A,x,\nabla u)\xi_{i}\xi_{j}\leq \Theta |\nabla u|^{p-2}_{A}|\xi|^{2}.$$
 On the other hand, for $p>2$ we apply the Cauchy-Schwarz inequality to (\ref{eq4_4}), and obtain
\begin{equation*}
\sum_{i,j=1}^{n}a^{ij}(A,x,\nabla u)\xi_{i}\xi_{j}\leq
|\nabla u|^{p-2}_{A}\Big[|\xi|^{2}_{A}+\frac{p-2}{|\nabla u|^{2}_{A}}|\nabla u|^{2}_{A}|\xi|^{2}_{A}\Big] =(p-1)\Theta|\nabla u|^{p-2}_{A}|\xi|^{2}_{A}.
\end{equation*}
\end{proof}
%%%%%%%%%%%%%
Proposition~\ref{bound} is crucial for the following regularity result. Indeed, since $A(x)$ satisfies \eqref{positive definite}, Lemma~\ref{NDP Lemma} below follows from \cite[Theorem~6.1.1]{PS}, \cite[Theorem~7.2, p.~290]{LU}, and \cite[Theorem 1]{GML}. We have
\begin{lemma}
\label{NDP Lemma}
Let $\Omega$ be a bounded domain in $\Real^{n}$ of class $C^{2,\alpha}$. Assume that the matrix $A$ is a bounded measurable symmetric matrix which is uniformly positive definite in $\Omega$, and $V\in L^\infty(\Gw)$.
Let $u\in W^{1,p}(\Omega)$ be a weak solution of the Dirichlet problem
\begin{equation} \label {NDP}{
\begin{cases}
Q_{A,V}(w)=f\geq 0 &\text{ in $\Omega$},\\
w=f_1\geq 0 &\text{  on $\partial\Omega$},
\end{cases}
}
\end{equation}
where $f\in L^{\infty}(\Omega) $, and $f_{1}\in C^{\alpha}(\partial \Omega)$ with $0<\alpha \leq 1$. Then $u\in L^\infty(\Omega)\cap C^\ga(\bar{\Omega})$.

If in addition the matrix $A\in C^{\gb}(\bar{\Gw})$ and $f_{1}\in C^{1,\gb}(\partial \Omega)$ for some $0<\gb\leq 1$, then there exists $0<\gg\leq 1$ such that $u\in C^{1,\gg}(\bar{\Omega})$.
\end{lemma}

%%%%%%%%%%%%%%%
Next, we generalize the results of J.~Garc\'ia-Meli\'an, and J.~Sabina de Lis \cite[Section 2]{GS} concerning the relationships between the principal eigenvalue, the weak and strong maximum principles, and the solvability of the Dirichlet problem.
\begin{definition}\label{def_SMP}{\em
By the {\em strong maximum principle} (SMP) of a quasilinear equation  $M(u)=0$ in $\Gw$  we mean the following statement: if $u$ is a nonnegative supersolution of the equation $M(u)=0$ in $\Gw$  with $u(x_{0})=0$ for some $x_{0} \in \Omega$, then $u= 0$ in $\Omega$.
 }
\end{definition}
Since the SMP is a local property, the weak Harnack inequality (see \cite[theorems~7.1.2 and 7.4.1]{PS}) clearly implies:
\begin{lemma}[SMP]\label{lem_SMP}
Assume that the matrix $A$ and the potential $V$ satisfy conditions (A), (E), and (V). If $\, 0\leq u\in W^{1,p}_\loc(\Gw)\cap C(\Gw)$  satisfies the differential inequality
\begin{equation}\label{eq_supers1}
    -\Delta_{p,A}(u)+V|u|^{p-2}u \geq 0 \qquad \mbox{in }\Gw,
\end{equation}
and $u(x_{0})=0$ for some $x_{0} \in \Omega$, then $u= 0$ in $\Omega$.
\end{lemma}
%%%%%%%%%%%%%%%%%%%%%%%
 We turn now to the boundary point lemma for our equation. First, we prove a boundary point lemma that holds for nonnegative functions satisfying the differential inequality
\begin{equation}\label{eq_dpa}
 -\Delta_{p,A}(u)+f(u)\geq 0 \qquad \mbox{in }\Gw,
\end{equation}
where, $f=f(u)$ satisfies the following condition (F):
\begin{equation}\tag{F}
 f \in C(0,\infty),\; f(0)=0, \;\mbox{and } f \mbox{ is nondecreasing on the interval } (0,\delta), \mbox{ where } \delta>0.
\end {equation}
\begin{theorem}[Boundary point lemma]\label{thm_SMP2}
Suppose that $A$ satisfies assumptions (A) and (E), and $f$ satisfy condition (F). Assume that $\Omega$ is of class $C^{2}$. Suppose that $f(s)>0$ for $s \in (0.\delta)$.
Let $x_{1}\in \partial \Omega$ satisfy the interior sphere condition, and assume that the matrix $A$ is of class $C^{2}$ in a closed relative neighborhood of $x_1$.
Let $u\in C^{1}(\bar{\Omega})$ be a positive solution of \eqref{eq_dpa} such that $u(x_1)=0$.   Then $\frac{\partial u(x_1)}{\partial \nu}<0$, where $\nu$ is the outer normal to $\partial\Omega$.
\end{theorem}
%%%%%%%%%%%%%%%%%%%%%%
\begin{proof} Let  $\mathcal{A}(s)=s^{p-2}$, and introduce the Riemannian metric $g^{ij}(x):=a_{ij}(x)$ induced by the matrix $A$, and let $\mathbf{s}(x):=\mathrm{dist}(x,x_0)$ be the induced geodesic distance from $x_0\in \Gw$. For certain ``radial" functions $v$ we need to estimate the expression $$\Delta_{p,A}(v)-f(v)=\partial_{x^{i}}\left\{g^{ij}(x)\mathcal{A}(|\nabla v|_{g})\partial_{x^{j}}v\right\}-f(v)$$
in the annular domain $G_S:=\{x\in\Gw\mid S/2< \mathbf{s}(x)< S\}$ centered at $x_0$.
 As in \cite[p.~220]{PS}, set $t:=S-\mathbf{s}(x)$, $\Phi(t):=t\mathcal{A}(t)$, and let $w$ be the (unique) solution of the problem
 \begin{equation*}
 \left\{
   \begin{array}{ll}
\displaystyle{-\dfrac{1}{(S-t)^k}\big[(S-t)^k\Gf(w'(t)\big]'+f(w(t))=0}& \qquad \mbox{in } (0,S/2),\\[4mm]
w(0)=0, w(S/2)=m>0,&
   \end{array}
 \right.
 \end{equation*}
where $k\in \mathbb{N}$ will be determined later.
The existence and uniqueness of $w$ is guaranteed by \cite[Lemma~4.2.1]{PS}. Moreover, by \cite[Lemma~4.2.2]{PS},  $w>0$ in $(0,S/2]$, and $w'>0$ in $[0,S/2]$.

 Consider the ``radial" function $v(x):=w(t)$.
By restricting the boundary value $v=m$ at $\partial{B_{S/2}}(x_0)$ to be sufficiently small, one can maintain $\sup{|\nabla v|_{g}} \leq \Theta|\nabla v| \leq 1$.

Using the summation convention, we have for the ``radial" function $v$ in $G_S$:
\begin{multline*}
\Delta_{p,A}(v)-f(v)=\partial_{x^{i}}\left\{g^{ij}(x)\mathcal{A}(|\nabla v|_{g})\partial_{x^{j}}v\right\}-f(v) =\\
g^{ij}(x)\partial_{x^{j}}\mathbf{s}(x)\partial_{x^{i}}\mathbf{s}(x)\left[\Phi(w')\right]'-\partial_{x^{i}}\left\{g^{ij}(x)\partial_{x^{j}}\mathbf{s}(x)\right\}\Phi(w')-f(w)  =\\
\left[\Phi(w')\right]'-\partial_{x^{i}}\left\{g^{ij}(x)\partial_{x^{j}}\mathbf{s}(x)\right\}\Phi(w')-f(w).
\end{multline*}
For $S$ small, we have  $\mathbf{s}\in C^2(G_S)$, and by \cite[Corollary~1.2]{SY}, there exists $\bar{k}\in \mathbb{N}$ such that
\begin{equation}\label{eq_Yau}
\Delta\mathbf{s}=\frac{1}{\sqrt{g(x)}}\partial_{x^{i}}\left\{\sqrt{g(x)}g^{ij}(x)\partial_{x^{j}}\mathbf{s}(x)\right\} \leq \frac{\bar{k}}{\mathbf{s}(x)}.
\end{equation}
A direct calculation shows that
\begin{multline*}
\partial_{x^{i}}\!\!\left\{g^{ij}(x)\partial_{x^{j}}\mathbf{s}(x)\!\right\}\!=\!
\frac{1}{\sqrt{g(x)}}\partial_{x^{i}}\!\!\left\{\sqrt{g(x)}g^{ij}(x)\partial_{x^{j}}\mathbf{s}(x)\!\right\}\!-\!
\frac{1}{\sqrt{g(x)}}\!\left\{g^{ij}(x)\partial_{x^{j}}\mathbf{s}(x)\!\right\}\!\partial_{x^{i}}\!\left\{\sqrt{g(x)}\!\right\}\!=\\
\Delta\mathbf{s}-\frac{1}{\sqrt{g(x)}}\left\{g^{ij}(x)\partial_{x^{j}}\mathbf{s}(x)\right\}\partial_{x^{i}}\left\{\sqrt{g(x)}\right\}.
\end{multline*}
Since $|\nabla \mathbf{s}(x)|$ is bounded, $g^{ij} \in C^{2}(\Omega)$, and $G_{S}$ is compact, we obtain using \eqref{eq_Yau} that there exists constant $k\in \mathbb{N}$ such that
\begin{equation}
\Delta_{p,A}\big(v(x)\big)-f\big(v(x)\big) \geq  \left[\Phi(w')\right]'-\frac{k}{\mathbf{s}}\Phi(w')-f(w)=0 \quad x\in G_S.
\end{equation}
We claim that $\Delta_{p,A}$ satisfies the monotonicity condition needed for the of the comparison principle \cite[Theorem~3.4.1]{PS}. Indeed,  $\boldsymbol{\mathcal{\hat{A}}(x,\boldsymbol{\xi})}:=g^{ij}(x)\mathcal{A}(|\boldsymbol{\xi}|_{g})\boldsymbol{\xi}$ satisfies \cite[(2.4.3)]{PS}:
\begin{multline}\label{mono2}
\left\langle \boldsymbol{\mathcal{\hat{A}}(x,\boldsymbol{\xi})}-\boldsymbol{\mathcal{\hat{A}}(x,\boldsymbol{\eta})},\boldsymbol{\xi}-\boldsymbol{\eta}\right\rangle= \\
\mathcal{A}(|\boldsymbol{\xi}|_{g})|\boldsymbol{\xi}|^{2}_{g}-\mathcal{A}(|\boldsymbol{\xi}|_{g})\left\langle [g^{ij}]\boldsymbol{\xi},\boldsymbol{\eta}\right\rangle-\mathcal{A}(|\boldsymbol{\eta}|_{g})\left\langle [g^{ij}]\boldsymbol{\eta},\boldsymbol{\xi}\right\rangle+\mathcal{A}(|\boldsymbol{\eta}|_{g})|\boldsymbol{\eta}|^{2}_{g} \geq \\
\left( \mathcal{A}(|\boldsymbol{\xi}|_{g})|\boldsymbol{\xi}|_{g}-\mathcal{A}(|\boldsymbol{\eta}|_{g})|\boldsymbol{\eta}_{g}|\right)
\left(|\boldsymbol{\xi}|_{g}-|\boldsymbol{\eta}|_{g}\right)=
 \left(|\boldsymbol{\xi}|^{p-1}_{g}-|\boldsymbol{\eta}|^{p-1}_{g}\right)\left(|\boldsymbol{\xi}|_{g}-|\boldsymbol{\eta}|_{g} \right) \geq 0,
\end{multline}
where we used the Cauchy-Schwarz inequality $|\!\left\langle [g^{ij}]\boldsymbol{\eta},\boldsymbol{\xi}\right\rangle\!| \leq |\boldsymbol{\eta}|_{g}|\boldsymbol{\xi}|_{g}$ and the monotonicity of $\Phi(t)=t^{p-1}$. Moreover, equality holds in \eqref{mono2} if and only if $\xi = \eta$.

Choosing $m$ sufficiently small, it follows from the comparison principle \cite[Theorem~3.4.1]{PS} that $u\geq v$ in $G_S$. Therefore, $\frac{\partial u}{\partial \nu}\leq\frac{\partial v}{\partial \nu}= -w'(t)<0$.

\end{proof}
Let   $M\geq \sup_{x\in \Gw}\left\{|V(x)|\right\}$. Then $f(u):=Mu^{p-1}$ satisfies condition (F), and  a nonnegative supersolution $u$ of the equation $Q_{A,V}(w)=0$ in $\Gw$ clearly satisfies
$$-\Delta_{p,A}(u) +Mu^{p-1}\geq 0.$$
%%%%%%%%%%%%%%%%%%%%%
Hence, Theorem~\ref{thm_SMP2} implies the following boundary point lemma for the operator $Q_{A,V}$:
\begin{Cor}[Boundary point lemma]\label{corSMP}
Suppose that the matrix $A$ and the potential $V$ satisfy assumptions (A), (E), and $(V)$.
Then the  following boundary point lemma  holds true (and also the SMP): Let  $u\gneqq 0$  satisfy the differential inequality
\begin{equation}\label{eq_supers}
    -\Delta_{p,A}(u)+V|u|^{p-2}u \geq 0 \qquad \mbox{in }\Gw.
\end{equation}
Then $u>0$ in $\Gw$. Suppose further that $u(x_1)=0$, where $x_{1}\in \partial \Omega$ satisfies the interior sphere condition,  and in a closed relative neighborhood $\Gw'$ of $x_1$:  $u\in C^{1}$, $V$ is bounded, and the matrix $A$ is uniformly positive definite and of class $C^{2}$.
Then $\frac{\partial u(x_1)}{\partial \nu}<0$, where $\nu$ is the outer normal to $\partial\Omega$.
\end{Cor}
%%%%%%%%%%%%%%%%
Next, we define a principal eigenvalue of the corresponding Dirichlet eigenvalue problem.
\begin{definition}{\em
For $\Gw\subset \R^n$, consider the eigenvalue problem
\begin{equation} \label{ev.pr}
\begin{cases}
Q_{A,V}(u)=\lambda |u|^{p-2}u &\text{ in $ \Gw$},\\
u\in W_{0}^{1,p}(\Gw), \; u\neq 0.&\\
\end{cases}
\end{equation}
We say that $\gl\in \R$ is a {\em principal eigenvalue} of  the operator $Q_{A,V}$ in $\Gw$ if  there exists  a nonnegative  function $u$ satisfying \eqref{ev.pr} (such $u$ is called a {\em principal eigenfunction}).
 }
\end{definition}
\begin{Rem} \label{L} {\em
Let $A\in C^2(\bar{\Gw},\R^{n^2})$ be a symmetric positive definite matrix, and let
$V \in L^{\infty}(\Omega)$, where $\Gw$ is a smooth bounded domain. If $\psi$ and $\phi$ are {\em nonnegative} eigenfunctions of problem (\ref{ev.pr}), then by the boundary point lemma (Corollary~\ref{corSMP}) we have $\frac{\psi}{\phi} \in L^{\infty}(\Gw)$.
 }
\end{Rem}

It turns out that if $\Gw\Subset \R^n$,  then  $Q_{A,V}$ admits a principal eigenvalue $\lambda_{1}$ defined by \eqref{p.ev} with a principal eigenfunction which is a minimizer of the variational problem:
\begin{equation} \label{p.ev}
\gl_1:=\gl_1(Q_{A,V},\Gw)= \gl_1(\Gw)=\inf_{u\in W_{0}^{1,p}(\Gw)}\frac{\int_{\Gw}(|\nabla u|_{A}^{p}+V(x)|u|^{p})\dx}{\int_{\Gw}|u|^{p}\dx}\,.
\end{equation}

\begin{lemma} \label{pos.e.f.}
Assume that the matrix $A\in L^\infty(\Gw,\R^{n^2})$ is a symmetric uniformly positive definite matrix in a bounded domain $\Omega$, and $V\in L^\infty(\Gw)$. Then the eigenvalue problem (\ref{ev.pr}) admits a principal eigenvalue $\gl$ with a principal eigenfunction $\phi\in W_{0}^{1,p}(\Gw)$. Such a principal eigenpair is given by
$\gl_1$, and a minimizer of \eqref{p.ev}. Furthermore, all eigenfunctions with eigenvalue $\lambda_{1}$ does not vanish in $\Gw$.

Moreover, if $A\in C^2(\bar{\Gw})$ and  $\Gw\in C^{2,\alpha}$, $0<\alpha<1$, then $\phi\in C^{1,\beta}\big(\overline{\Gw}\big)$  for some $0<\beta<1$, and  $\,\frac{\partial \phi}{\partial \nu}<0$ on $\partial \Gw$.
\end{lemma}

\begin{proof}
We repeat, with necessary changes, the proof of \cite[Lemma 3]{GS}.
Since $V\in L^{\infty}(\Gw)$, we may assume that $V(x)\geq0$. Otherwise, replace $V$ with $V_{M}(x):=V(x)+M\geq 0$ and $\lambda_{M}:=\lambda+M$.

Now, if $V\geq 0$, then the functional $\mathcal{Q}_{A,V}$
is sequentially weakly lower semicontinuous in $W_{0}^{1,p}(\Gw)$, and coercive in
$\mathcal{V}:=\{u\in W_{0}^{1,p}(\Gw) | \int_{\Gw}{|u|^{p}dx}=1\}.$
Hence, the infimum $\lambda_{1}$ in \eqref{p.ev} is attained.
In particular, there exists $\phi \in \mathcal{V}$ such that $\phi$ is a weak solution of the equation
\begin{equation} \label{ev.eq}
-\Delta_{p,A}(u)+V(x)|u|^{p-2}u=\lambda_{1} |u|^{p-2}u.
\end{equation}
Since $|\nabla (|\phi|)|\leq|\nabla \phi|$, we get that $|\phi|$ is also a minimizer of (\ref{p.ev}) and hence it is a nonnegative weak solution of equation \eqref{ev.eq}.
Thus, by the Harnack inequality either $|\phi|>0$ or $|\phi|=0$, and consequently, $\phi$ does not vanish in $\Gw$, and has a definite sign in $\Gw$. Furthermore, by the same reasoning all eigenfunctions with eigenvalue $\lambda_{1}$ does not vanish in $\Gw$.

Moreover,  if $A\in C^{\ga}(\bar{\Gw})$, then by Lemma~\ref{NDP Lemma}, $\phi \in C^{1,\beta}\big(\overline{\Gw}\big)$.
Finally, if $A\in C^2(\bar{\Gw})$, then by the SMP and the boundary point lemma (Corollary~\ref{corSMP}), we have that $\frac{\partial |\phi|}{\partial \nu}<0$ on $\partial \Gw$.
\end{proof}

Using our previous results we extend the main theorem of J. Garc\'ia-Meli\'an, and J. Sabina de Lis \cite[Theorem~2]{GS}.   We have:
\begin{theorem}
\label{Weak, Strong, Pos.}
Let $A$ be a bounded symmetric matrix which is uniformly positive definite in a $C^{1,\ga}$-bounded domain $\Omega$ ($0<\alpha\leq 1$), and $V\in L^\infty(\Gw)$.  Consider the following assertions.
\begin{itemize}
 \item[(i)] $Q_{A,V}(u)$ satisfies the following weak maximum principle (WMP): If $u \in W^{1,p}(\Omega)$ is a solution of the equation $Q_{A,V}(u)=f\geq0$ in $\Omega$ with some $f\in L^{\infty}(\Omega)$, and satisfies $u\geq0$ on $\partial \Omega$, then $u\geq0$ in $\Omega$.
 \item[(ii)] $Q_{A,V}(u)$ satisfies the following version of the strong maximum principle: If $u \in W^{1,p}(\Omega)$ is a solution of the equation $Q_{A,V}(u)=f\geq 0$ in $\Omega$ with some $f\in L^{\infty}(\Omega)$, and satisfies $u\geq0$ on $\partial \Omega$, and $u\neq0$ in $\Omega$, then $u>0$ in $\Omega$.
 \item[(iii)] $\lambda_1(\Omega)>0$, where $\gl_1$ is defined by \eqref{p.ev}.
 \item[(iv)] There exists a positive strict supersolution $v \in W^{1,p}_{0}(\Omega)$ of $Q_{A,V}(w)=0$ in $\Gw$, with $Q_{A,V}(v) \in L^{\infty}(\Omega)$, i.e., $Q_{A,V}(v)=f$ in $\Omega$, where $f\in L^{\infty}(\Omega), f\gneqq0$.
 \item[(iv')] There exists a positive strict supersolution $v \in W^{1,p}(\Omega)\cap L^{\infty}(\Omega)$  of $Q_{A,V}(w)=0$ in $\Gw$, with  $Q_{A,V}(v) \in L^{\infty}(\Omega)$, i.e., $Q_{A,V}(v)=f$ in $\Omega$, where $f\in L^{\infty}(\Omega), f\gneqq0$.
\item[(v)] For each nonnegative $f\in L^{\infty}(\Omega)$ there exists a nonnegative weak solution $u \in W^{1,p}_{0}(\Omega)$ of the problem $Q_{A,V}(w)=f$ in $\Omega$, and $w=0$ on $\partial \Omega$.
\end{itemize}
Then $(i)\Leftrightarrow (ii)\Leftrightarrow (iii)\Rightarrow(iv)\Rightarrow (iv')$,  and  $(iii)\Rightarrow (v)\Rightarrow (iv')$.
\end{theorem}
\begin{proof}
$(i)\Longrightarrow (ii).$ Let $u \in W^{1,p}(\Omega)$ be a solution of the equation $Q_{A,V}(w)=f \geq 0$ in $\Omega$ with $f\in L^{\infty}(\Omega)$ such that $u\geq0$ on $\partial \Omega$, and $u\neq 0$  in $\Omega$. The WMP {\em (i)} implies that  $u \geq 0$ in $\Omega$. Hence, the strong maximum principle of Lemma~\ref{lem_SMP} implies that $u>0$ in $\Omega$.

$(ii)\Longrightarrow (iii).$ Assume that $\gl_1 \leq 0$, and let $\phi>0$ be the corresponding principal eigenfunction. Then $\psi:=-\phi$ satisfies $$-\Delta_{p,A}(\psi)+V|\psi|^{p-2}\psi=\gl_1 |\psi|^{p-2}\psi \geq 0,$$ and $\psi=-\phi=0$ on $\partial\Omega$. By {\em (ii)} we have $\psi>0$ in $\Omega$, hence, $\phi<0$, a contradiction.

$(iii)\Longrightarrow (i).$ Assume, that $u \in W^{1,p}(\Omega)$ is a solution of the equation $Q_{A,V}(w)=f\geq0$ in $\Omega$ and $u\geq 0$ on $\partial \Omega$, with $f\in W^{-1,p'}(\Omega)$, and  $\left\langle f,v \right\rangle \geq 0$ for each nonnegative $v \in W_{0}^{1,p}(\Omega)$.
Denote $u_{-}(x):=\min\{u(x),0\}$, hence $u_{-} \in W_{0}^{1,p}(\Omega)$. Consequently,
\begin{multline*}
\int_{\Omega}(|\nabla u_{-}|_{A}^{p}+V|u_{-}|^{p})\dx=
\int_{\Omega}(|\nabla u_{-}|_{A}^{p-2}A(x)\nabla u_{-} \cdot \nabla u_{-}+V|u_{-}|^{p-2}u_{-}u_{-})\dx=\\[2mm]
\int_{\Omega}(|\nabla u|_{A}^{p-2}A(x)\nabla u \cdot \nabla u_{-}+V|u|^{p-2}uu_{-})\dx=
 \left\langle f,u_{-} \right\rangle \leq 0< \lambda_{1}.
\end{multline*}
In light of the definition of $\lambda_{1}$, we obtain that $u_{-}= 0$, so, $u\geq 0$ in $\Omega.$

$(iii)\Longrightarrow (iv).$ The principal eigenfunction is a desired positive supersolution.

$(iv)\Longrightarrow (iv').$ This implication is trivial.

$(iii)\Longrightarrow (v).$ Let $f\in W^{-1,p'}(\Omega)$, and define a functional
$$\mathcal{J}_{f}(u)=\int_{\Omega}(|\nabla u|_{A}^{p}+V|u|^{p}-pfu)\dx.$$
Since $\lambda_{1}>0$, the functional $\mathcal{J}_{f}(u)$ is coercive in $W_{0}^{1,p}(\Omega)$.
Hence, for each $0\leq f \in L^{\infty}(\Omega)$ there exists a weak solution $u \in W_{0}^{1,p}(\Omega)$ of problem (\ref{NDP}) (with $f_1=0$). Moreover, since $(iii)\Longrightarrow (i)$ and $(ii)$, it follows that $u\geq 0$ in $\Gw$, and if $f \neq 0$ it follows that $u>0$ in $\Gw$.

$(v)\Longrightarrow (iv').$ Follows directly from the SMP.
\end{proof}
\begin{remark}
 {\em Later (in Corollary~\ref{cor_simple1}) we show that $(iv')\Rightarrow (iii)$. Hence, all the assertions of Theorem~\ref{Weak, Strong, Pos.} are in fact equivalent.}
\end{remark}
%%%%%%%%%%%%%%%%%%%%%%
For the Dirichlet problem \eqref{NDP} we have the following uniqueness result.
%%%%%%%%%%%%%%%%%%%%%%%%%
\begin{theorem}
\label{Weak, Strong, Pos.S}
Assume that $\Omega \subset \Real^{n}$ is a bounded $C^{2,\alpha}$-domain, $0<\alpha<1$, $A\in C^{2}(\bar{\Gw},\R^{n^2})$ is positive definite, and $V\in L^{\infty}(\Omega)$. Then for nonnegative $f\in L^{\infty}(\Omega)$ and $f_1 \in C^{1,\ga}(\partial \Omega)$ there exists at most one nonzero nonnegative weak solution $u \in W^{1,p}(\Omega)\cap L^\infty(\Omega)$ of problem \eqref{NDP} (up to a multiplicative constant if $f=f_1=0$).

Moreover, if $f_1>0$, $f_1 \in C^\ga(\partial \Omega)$, then the uniqueness of \eqref{NDP} holds even if $A$ is only a bounded measurable symmetric matrix which is uniformly positive definite.
\end{theorem}
\begin{proof}
Let $u,v \in W^{1,p}(\Omega)\cap L^\infty(\bar{\Omega})$ be nonnegative solutions of (\ref{NDP}).  In light of
\cite[Theorem~7.2, p.~290]{LU}, $u,v \in C^\ga(\bar{\Gw})$. By the SMP (Lemma~\ref{lem_SMP}), any such solution is positive in $\Gw$. Hence, (using the boundary point lemma, in case $f_1 \not>0$ and $A\in C^2$), we have  $\frac{v}{u}, \frac{u}{v} \in L^\infty(\Omega)$.
Therefore, $(u,v) \in \mathcal{D}(I)$, and by Lemma~\ref{I(u,v)} we have
\begin{multline*}
0 \leq \left\langle -\Delta_{p,A}(u),\frac{u^p-v^p}{u^{p-1}}\right\rangle_{\Omega}-\left\langle -\Delta_{p,A}(v),\frac{u^p-v^p}{v^{p-1}}\right\rangle_{\Omega}=\\
   \left\langle\!\! f-V|u|^{p-2}u,\frac{u^p-v^p}{u^{p-1}}\!\right\rangle_{\Omega}\!\!-\left\langle\!\! f-V|v|^{p-2}v,\frac{u^p-v^p}{v^{p-1}}\!\right\rangle_{\Omega}\!=\!
   \int_{\Omega}\!\!f\!\frac{(u^p-v^p)(v^{p-1}-u^{p-1})}{u^{p-1}v^{p-1}}\dx \leq 0. \nonumber
\end{multline*}
It follows that $I(u,v)=0$, and Lemma~\ref{I(u,v)} implies that $u=kv$ for some positive $k$,
If either $f \gneqq 0$ or $f_1\gneqq 0$, it follows that $k=1$, and hence $u=v$.
\end{proof}
The following is an extension of the weak comparison principle (WCP) \cite[Theorem 5]{GS}.
\begin{theorem} \label{compare thm}
Assume that $\Omega \subset \Real^{n}$ is a bounded $C^{2,\alpha}$-domain, $0<\alpha<1$, $A$ is a symmetric bounded matrix which is uniformly positive definite in $\Gw$, and $V\in L^{\infty}(\Omega)$.
Assume that $\lambda_1>0$, where $\gl_1$ is defined by \eqref{p.ev}. Let $u_{i}\in W^{1,p}(\Omega)\cap C^\ga(\bar{\Omega})$ satisfy $Q_{A,V}(u_{i})\in L^{\infty}(\Omega)$, where $i=1,2$. Suppose further that the following inequalities are satisfied
\[
\begin{cases}
Q_{A,V}(u_{1})\leq Q_{A,V}(u_{2}) &\text{in $\Omega$},\\
Q_{A,V}(u_{2})\geq0 &\text{in $\Omega$},\\
u_{1}\leq u_{2} &\text{on $\partial\Omega$},\\
u_{2}> 0 &\text{on $\partial\Omega$}.
\end{cases}
\]
Then, $u_{1}\leq u_{2}$ in $\Omega$.
Moreover, if the conditions for the validity of the boundary point lemma are satisfied, then the conclusion holds true even if $u_2\geq 0$ on $\partial\Omega$.
\end{theorem}
\begin{proof}
Since $u_{2}>0$ in $\bar{\Omega}$, there exists a constant $c>1$ such that $u_{1}<cu_{2}$ in $\Omega$. Set $g:=Q_{A,V}(u_{2})$ , $g_{2}=u_{2}|_{\partial \Omega}$, and consider the Dirichlet problem
\begin{equation} \label{sub-super} {
\begin{cases}
Q_{A,V}(w)=g &\text{ in $\Omega$},\\
w=g_{2} &\text{ on  $\partial\Omega$}.
\end{cases}
}
\end{equation}
Clearly, $u_{1}$ is a subsolution and $cu_{2}$ is a supersolution of (\ref{sub-super}). Therefore, by the sub/super-solution technique (see \cite[Theorem~4.14, p.~272]{DJ}), there exists a weak solution $v \in W^{1,p}(\Omega) \cap L^{\infty}(\Omega)$ to (\ref{sub-super}), satisfying $u_{1} \leq v \leq cu_{2}$, and $v=g_{2}|_{\partial \Omega}$ in the sense of traces. Moreover, Theorem~\ref{Weak, Strong, Pos.} implies that $v > 0$ in $\Omega$.  In light of Lemma~\ref{NDP Lemma}, $v\in C^\ga(\bar{\Gw})$ (or even $v\in C^{1,\ga}(\bar{\Gw})$ if $A\in C^\ga(\bar{\Gw})$). Consequently, by uniqueness (Theorem~\ref{Weak, Strong, Pos.S}), $v=u_{2}$. Thus, $u_{1}\leq u_{2}$.
\end {proof}
\begin{problem}
{\em Prove the WCP assuming only $v_2\geq 0$ on $\partial\Gw$ (without using the boundary point lemma).
 }
\end{problem}

%%%%%%%%%%%%%%%%%%%%%%%%%%%%%%%
\mysection{The Agmon-Allegretto-Piepenbrink (AAP) theory}\label{sec_AAP}
In the present section we generalize the AAP theorem claiming that $\mathcal{Q}_{A,V}\geq 0$ in $\Gw$ if and only if
the equation $Q_{A,V}(u)=0$ admits a positive (super)solution in $\Gw$.
First, we need to prove the strict monotonicity of  $\lambda_{1}$ with respect to the domain (see \cite[Theorem~2.3]{AH1}).
\begin{lemma} \label{mono}
Let $\Omega_{1} \Subset \Omega_{2}\Subset \Gw$ be smooth bounded domains, and suppose that $A$ and $V$ satisfy assumptions (A), (E), and (V) in $\Gw$, and $\mathcal{Q}_{A,V} \geq 0$ in $\Omega_{2}$. Then $\lambda_{1}(\Omega_{1}) > \lambda_{1}(\Omega_{2}) \geq 0$.
\end{lemma}
%%%%%%%%%%%%%%%%%%%%%%
\begin{proof}
It follows from \eqref{p.ev} that $\lambda_{1}(\Omega_{2}) \geq 0$. Let $\gf_i\in W^{1,p}_{0}(\Omega_{i})$ be normalized principal eigenfunctions of the operator $Q_{A,V}$ in $\Gw_i$ with eigenvalues $\lambda_{1}(\Omega_{i})$, $i=1,2$. Let $\{\vgf_k\}\subset C_0^\infty(\Gw_1)$ be a nonnegative minimizing sequence that converges to $\phi_1$ in  $W^{1,p}_{0}(\Omega_{1})$. By Picone identity we have
\begin{multline*}\label{eq_strict_mon}
    0\leq \int_{\Gw_2}L_A(\vgf_k,\gf_2)\dx= \int_{\Gw_1}L_A(\vgf_k,\gf_2)\dx= \int_{\Gw_1}R_A(\vgf_k,\gf_2)\dx=\\
    \int_{\Gw_1}\!\!\!\!|\nabla \vgf_k|_{A}^{p}\!\dx-\!\int_{\Gw_1}\!\!\!\!\nabla \!\left(\frac{\vgf_k^{p}}{\gf_2^{p\!-\!1}}\right)\cdot A(x)\nabla \gf_2|\nabla \gf_2|_{A}^{p\!-\!2}\!\dx\!=\!\!
    \int_{\Gw_1}\!\!\!\!|\nabla \vgf_k|_{A}^{p}\!\dx+\!\int_{\Gw_1}\!\!\!\! V(x)\vgf_k^{p}\!\dx-\!\lambda_{1}(\Omega_{2}) \!\!\int_{\Gw_1} \!\!\!\!\vgf_k^{p}\!\dx.
\end{multline*}
Letting $k\to \infty$, and using the Fatou's lemma,  we arrive at
\begin{equation*}
0\leq \int_{\Gw_1}L_A(\gf_1,\gf_2)\dx\leq \lim_{k\to \infty} \int_{\Gw_1}L_A(\vgf_k,\gf_2)\dx =\Big(\lambda_{1}(\Omega_{1})-\lambda_{1}(\Omega_{2})\Big)  \int_{\Gw_1} \gf_1^p\dx.
\end{equation*}
So, $\lambda_{1}(\Omega_{1})\geq \lambda_{1}(\Omega_{2})$.
If $\lambda_{1}(\Omega_{1})=\lambda_{1}(\Omega_{1})$, then $L_A(\gf_1,\gf_2)=0$ a.~e. in $\Gw_1$. Hence, Proposition~\ref{Picone} implies that there is a constant $c> 0$ such that $\gf_2\mid_{\Gw_1}=c\gf_1$, and this is impossible since $\gf_2>0$ on $\partial \Gw_1$ (by the Harnack inequality). Thus, $\lambda_{1}(\Omega_{1})>\lambda_{1}(\Omega_{1})$.
\end{proof}
Using our earlier results, we extend now the AAP-type theorem.
\begin{theorem}[AAP-type theorem]\label{pos sol}
Suppose that $A$ and $V$ satisfy assumptions (A), (E), and (V) in a domain $\Gw\subset \R^n$. Then the following assertions are equivalent:

\begin {itemize}
 \item [(i)] The functional $\mathcal{Q}_{A,V}$ is nonnegative on $C_{0}^{\infty}(\Omega)$.
 \item [(ii)] The equation $Q_{A,V}(u)=0$ in $\Gw$ admits a positive solution.
 \item [(iii)] The equation $Q_{A,V}(u)=0$ in $\Gw$  admits a positive supersolution.
\end{itemize}

\end{theorem}

\begin{proof}
$(i)\Longrightarrow (ii)$. Assume that $\mathcal{Q}_{A,V} \geq 0$ on $C_{0}^{\infty}(\Omega)$, and let $\{\Omega_{N}\}$ be an exhaustion of $\Gw$.
By Lemma~\ref{mono}, $\lambda_{1}(\Omega_{N})>0$ for all $N\in \mathbb{N}$. Let $f_{N} \in C^{\infty}_{0}(\Omega_{N}\setminus \Omega_{N-1})$ be a nonnegative nonzero function.
Theorem~\ref{Weak, Strong, Pos.} implies the existence of a positive solution $v_{N}$ of the problem
\[
\begin{cases}
Q_{A,V}(w)=f_{N} &\text{in $\Omega_{N}$},\\
w=0 &\text{on $\Omega_{N}$}.
\end{cases}
\]
Set $u_{N}(x):=\frac{v_{N}(x)}{v_{N}(x_{0})}$. By the Harnack convergence principle, $\{u_{N}\}$ admits a subsequence which converges locally uniformly to a positive solution $u$ of the equation $Q_{A,V}(w)=0$ in $\Omega$.

$(ii)\Longrightarrow (iii)$. This implication is trivial.

$(iii)\Longrightarrow (i)$. Suppose that $v$ is a positive supersolution of (\ref {derivative}) in $\Gw$. Then Proposition~\ref{Q=L} implies that $\mathcal{Q}_{A,V}(\vgf) \geq 0$ for all {\em nonnegative} $\vgf \in C^{\infty}_{0}(\Omega)$. Since $\mathcal{Q}_{A,V}(u)=\mathcal{Q}_{A,V}(|u|)$ on $C^{\infty}_{0}(\Omega)$, a standard approximation argument shows that  $\mathcal{Q}_{A,V} \geq 0$ on $C^{\infty}_{0}(\Omega)$.
\end{proof}
%%%%%%%%%%%%%%%%%%%%%%%%%%%%%%%%%%%%%%%%%%%%%%%%%%%%%%%
\mysection{The Main Theorem}\label{sec_main_thm}
This section is devoted to the following result which generalizes \cite[Theorem~3.3]{PT2}.
\begin{theorem}[Main Theorem] \label{main thm}
Suppose that the matrix $A$ and the potential $V$ satisfy conditions (A), (E) and (V). If $p<2$ assume further that $A\in C^\ga(\Gw)$. Suppose that $\mathcal{Q}_{A,V}$ is nonnegative on $C^{\infty}_{0}(\Omega)$. Then
\begin{enumerate}
	\item [1.] Any ground state $\phi$ is a positive solution of (\ref{derivative}).
	\item [2.] $\mathcal{Q}_{A,V}$ is critical in $\Omega$ if and only if $\mathcal{Q}_{A,V}$ admits a null sequence. Moreover, there exists a null sequence that converges locally uniformly in $\Omega$ to the ground state.
	\item [3.] $\mathcal{Q}_{A,V}$ admits a null sequence if and only if (\ref{derivative}) admits a unique positive continuous supersolution. \item [4.] $\mathcal{Q}_{A,V}$ is subcritical in $\Gw$ if and only if there exists a strictly positive continuous function $W$ such that $\mathcal{Q}_{A,V-W}$ is nonnegative on $C^{\infty}_{0}(\Omega)$.
	\item [5.] If $\mathcal{Q}_{A,V}$ admits a ground state $\phi$, then the
following Poincar\'{e} type inequality holds:

There exists $0<W\in C(\Gw)$ such that for every $\psi\in C_{0}^{\infty}(\Omega)$
satisfying $\int_\Gw \psi \phi\dx\neq 0$ there exists a constant $C>0$ so that
the following inequality holds:
\begin{equation} \label{poincare}
\mathcal{Q}_{A,V}(\vgf)+C\Bigl|\int_\Gw \vgf\,\psi \dx\Bigr|^{p}\geq C^{-1}\int_{\Omega}W(x)|\vgf|^{p}\dx\quad\forall \vgf\in C_{0}^{\infty}(\Omega).
\end{equation}
\end{enumerate}
\end{theorem}

The proof of Theorem~\ref{main thm} is based on the proofs in \cite[Section~3]{PT1} and \cite[Theorem~4.3]{PT3}. First we need to generalize Lemma~3.1 and Lemma~3.2 in \cite{PT1}.

Let $B \subset \Omega$ be a nonempty open set, and set
\begin {equation}
c_{B}:=\inf_{\substack{\\\vgf \in C_{0}^{\infty}(\Omega)\\[1mm] \int_{B}|\vgf|^{p}\dx=1}}\mathcal{Q}_{A,V}(\vgf)= \inf_{\substack{\\0 \leq \vgf \in C_{0}^{\infty}(\Omega)\\[1mm] \int_{B}|\vgf|^{p}\dx=1}}\mathcal{Q}_{A,V}(\vgf).
\end{equation}
Clearly, the criticality of $\mathcal{Q}_{A,V}$ in $\Gw$ implies that $c_B=0$ for any nonempty open set in $\Gw$.
\begin{lemma} \label{subcrit}
If for every open set $B \Subset \Omega, c_{B}>0 $, then there exists a continuous positive function $W$, such that
\begin{equation} \label{W1}
\mathcal{Q}_{A,V}(\varphi) \geq \int_{\Omega}{W(x)|\varphi(x)|^{p}\dx} \qquad \forall \varphi \in C_{0}^{\infty} (\Omega).
\end{equation}
\end{lemma}

\begin{proof} The proof is obtained by a partition of unity argument as in \cite[Lemma~3.1]{PT1}.
\end{proof}

\begin{lemma} \label{crit}
Suppose that $\mathcal{Q}_{A,V}\geq 0$ in $\Omega$, where $A$ and $V$ satisfy conditions (A), (E) and (V). If $p<2$ assume further that $A\in C^\ga(\Gw)$.
If there exists a nonempty open set $B \Subset \Omega$ such that $c_{B}=0$, then $\mathcal{Q}_{A,V}$ admits a ground state $\phi$. Moreover, $\phi$ is the unique positive supersolution of the equation $Q_{A,V}(u)=0$ in $\Gw$.
\end{lemma}
%%%%%%%%%%%%%%%%
\begin{proof}
Since $c_{B}=0$, it follows that there exists a null sequence $\{\vgf_{k}\}\subset C_{0}^{\infty}(\Omega)$,  such that $\vgf_k\geq 0$, and $\int_{B}{|\vgf_{k}|^p\dx}=1$. Fix a positive (super)solution $v\in W^{1,p}_\mathrm{loc}(\Gw)$ of \eqref{derivative}.

\medskip

Assume first that $p\geq 2$. So, we may assume that $v\in C^\ga(\Gw)$. Denote $w_k:=\vgf_k/v$. By Lemma~\ref{lem_simplified} we have
$\int_{\Omega}v^{p}|\nabla w_k|^{p}_{A}\dx\leq C \mathcal{Q}_{A,V}(\vgf_k) \to 0$. Hence,  $\nabla (w_k)\to 0$ in $L^p_\mathrm{loc}(\Gw)$.
Poincar\'{e} inequality in $C^1$-subdomains  $\gw\Subset \Gw$, implies that $w_k \to \mathrm{const}$ in $W^{1,p}_\mathrm{loc}(\Gw)$.  Consequently, there exists $c\geq 0$ such that (up to a subsequence) $\vgf_k\to cv$ a.e. in $\Gw$, and also in $L^p_\mathrm{loc}(\Gw)$. Thus,  $c^{-p}=\int_B v^p\dx>0$. It follows that any null sequence $\{\vgf_k\}$ converges to the same positive supersolution $v$. Hence, $v$ is the unique positive solution, and also the unique positive supersolution of the equation $Q_{A,V}(u)=0$ in $\Gw$.

\medskip

Assume now that $1<p< 2$. In this case the proof is more involved and consists of three steps (see the proof of \cite[Lemma~3.2]{PT1}).

\noindent \textbf{Step 1}.
By our assumption,  we may assume that $v\in C^{1,\ga}(\Gw)$. Let $\omega \Subset \Omega$ be an open connected set containing $B$, and let $\omega' \subset \omega$. By \eqref{L_A} we have that
\begin{equation*}
\int_{\omega'}\!\!\!L_{A}(\vgf_k,v)\dx\! = \!\!\int_{\omega'}\!\!|\nabla \vgf_k|^{p}_{A}\dx+(p\!-\!1)\!\!\int_{\omega'}\!\!\left(\frac{\vgf_k}{v}\right)^{p}\!|\nabla v|^{p}_{A}\dx -
p\!\int_{\omega'}\!\!\left(\frac{\vgf_k}{v}\right)^{p\!-\!1}\!\nabla \vgf_k \cdot A(x)\nabla v|\nabla v|^{p\!-\!2}_{A}\dx.
\end{equation*}
Recall that $L_{A}(\vgf_k,v)\! \geq\! 0$. Thus, by the Cauchy-Schwarz inequality we obtain
$$\int_{\omega'}|\nabla \vgf_k|^{p}_{A}\dx \leq \int_{\Omega}L_{A}(\vgf_k,v)\dx+p\int_{\omega'}|\nabla \vgf_k|_{A}\left(\frac{\vgf_k}{v}|\nabla v|_{A}\right)^{p-1}\!\!\dx -(p\!-\!1)\!\!\int_{\omega'}\!\!\left(\frac{\vgf_k}{v}\right)^{p}\!|\nabla v|^{p}_{A}\dx.$$
Young's inequality, $ab\leq \dfrac{a^{p}}{p}+\dfrac{b^{p'}}{p\ \!'}$,  with $a:=(1\!-\!\vge)^{1/p}|\nabla \vgf_k|_{A}, b:= (1\!-\!\vge)^{-1/p}\left(\frac{\vgf_k}{v}|\nabla v|_{A}\right)^{p-1}$ implies
$$\vge\int_{\omega'}|\nabla \vgf_k|^{p}_{A}\dx \leq
\int_{\Omega}L_{A}(\vgf_k,v)\dx+(p\!-\!1)\left(\frac{1}{(1-\vge)^{1/(p-1)}} -1 \right) \int_{\omega'}\!\!|\nabla(\log  v)|_{A}^{p}\vgf_k^{p}\dx.$$
By Proposition~\ref{Q=L},  $\int_{\Omega}L_{A}(\vgf_k,v)\dx  = o(1)$, where $o(1)\to 0$ as $k\to \infty$. Hence,
\begin{equation} \label{50}
\int_{\omega'}|\nabla \vgf_k|^{p}_{A}\dx \leq o(1)+C(\omega)\int_{\omega'}\vgf^{p}_{k}\dx.
\end{equation}
%%%%%%%%%%%%%%%%%%%%%%%%%%%%%%%
\begin{equation*}
\mbox{\textbf{Step 2}. Set:}\qquad  \omega_{0}:=\left\{x\in \omega \mid \exists\, \rho(x)\in\left(0,d(x,\Omega \backslash \omega)\right),\ \sup_{k}\int_{B_{\rho(x)}(x)}|\vgf_k|^{p}\dx < \infty\right\}.
\end{equation*}
Since $\int_{B}|\vgf_k|^{p}\dx=1$, we get that $B \subset \omega_{0}$.
Moreover, $\omega_{0}$ is readily an open set. We claim that $\omega_{0}$ is relatively closed  in $\omega$. Indeed,
we use the following version of Poincar\'e inequality \cite[Theorem~8.11]{LL}. For fix $0<r<1$ we have
$$\int_{B_{1}(0)}\left|u-\dfrac{\int_{B_{r}(0)}u\dy}{\gw_nr^n}\right|^{p}\dx \leq  C(p,n,r)\int_{B_{1}(0)}|\nabla u|^{p}\dx \qquad \forall u\in W_\loc^{1,p}(\R^n).$$
On the other hand, there exists  $\theta:=\theta_B>0$ such that $|\nabla u|^{p} \leq \theta^{-\frac{p}{2}}|\nabla u|^{p}_{A}$. Hence, in light of the inequality
$2^{1-p}|a|^{p}-|b|^{p}\leq |a-b|^{p}$, we get
\begin{equation*}  2^{1-p}\int_{B_{1}(0)}\!|u|^{p}\dx-\int_{B_{1}(0)}\left|\dfrac{\int_{B_{r}(0)}u\dy}{\gw_nr^n}\right|^{p}\dx \leq \theta^{-\frac{p}{2}}C(p,n,r)\int_{B_{1}(0)}\!|\nabla u|^{p}_{A}\dx.
\end{equation*}
Hence, for all $ u \in W^{1,p}(\R^n)$ we obtain
\begin{equation} \label{51}
\int_{B_{1}(0)}|u|^{p}\dx  \leq  \tilde{C}\int_{B_{1}(0)}|\nabla u|^{p}\dx + \frac{C}{r^{n}}\int_{B_{r}(0)}|u|^{p}\dx.
\end{equation}
By scaling, we obtain for $0<a<1$ that
\begin{equation*}
\int_{B_{a}(0)}|u|^{p}\dx \leq \\ \tilde{C}a^{p}\int_{B_{a}(0)}|\nabla u|^{p}\dx + Cr^{-n}\int_{B_{ar}(0)}|u|^{p}\dx \qquad \forall u\in W_\loc^{1,p}(\R^n).
\end{equation*}
Hence, for every $\varepsilon > 0$ and $\rho > 0$ small enough, take $a_{\varepsilon}<\min\left\{\left(\frac{\varepsilon}{\tilde{C}}\right)^{1/p},\frac{\rho}{r}\right\}$, so that
\begin{equation}\label{52}
\int_{B_{a_{\varepsilon}}(0)}|u|^{p}\dx \leq \varepsilon\int_{B_{a_{\varepsilon}}(0)}|\nabla u|^{p}_{A}\dx+C(\varepsilon,\rho)\int_{B_{\rho}(0)}|u|^{p}\dx.
\end{equation}
Therefore, we get for every $x \in \omega$, $\vge>0$, $\delta \in \left(0,\min\left\{a_{\varepsilon},d(x,\partial{\omega})\right\}\right)$, and $u \in W^{1,p}_{0}(\Gw)$, that the following inequality holds
\begin{equation} \label{53}
\int_{B_{\delta}(x)}|u|^{p}\dx \leq \varepsilon\int_{B_{\delta}(x)}|\nabla u|^{p}_{A}\dx+C(\varepsilon,\rho)\int_{B_{\rho}(x)}|u|^{p}\dx.
\end{equation}

Let $x_{j} \in \omega_{0},\ x_{j}\rightarrow x_{0} \in \omega$. Let $\varepsilon < \frac{1}{2C(\omega)}$, where $C(\omega)$ is the constant in (\ref{50}). Let $a_{\varepsilon}$ be as in (\ref{52}) and fix $\delta \in \left(0,\min\left\{a_\vge,d(x,\partial{\omega})\right\}\right)$, and $u \in W^{1,p}_{0}(\Omega)$. Finally, fix $j$ such that $|x_{0}-x_{j}|<\frac{\delta}{2}$. Then, with $\rho=\rho(x_{j})$, we get from (\ref{50})  and (\ref{53}) that
\begin{equation} \label{54}
\frac{1}{2C(\omega)}\int_{B_{\delta}(x_{j})}|\nabla \vgf_k|^{p}_{A}\dx \leq C(\varepsilon,\rho(x_{j}))\int_{B_{\rho(x_{j})}(x_{j})}|\vgf_k|^{p}\dx+o(1),
\end{equation}
where $o(1)\to 0$ as $k\to\infty$. Thus $\int_{B_{\delta}(x_{j})}|\nabla \vgf_k|^{p}_{A}\dx$ is bounded (in $k$), and by (\ref{53}), $\int_{B_{\delta}(x_{j})}|\vgf_k|^{p}\dx$ is also bounded. By the choice of $j$, we have that $B_{\frac{\delta}{2}}(x_{0}) \subset B_{\delta}(x_{j})$, it follows that $\int_{B_{\frac{\delta}{2}}(x_{0})}|\vgf_k|^{p}\dx$ is bounded and consequently, $x_{0} \in \omega_{0}$.So, $\omega_{0}$ is also relatively closed in $\omega$.
Since $\omega$ is connected, we obtain  $\omega_{0}=\omega$. Hence, $\left\{\vgf_k\right\}$ is bounded in $L^{p}_{\loc}(\Omega)$, and by (\ref{50}) it follows that $\left\{\vgf_k\right\}$ is bounded in $W^{1,p}_{\loc}(\Omega)$.

\noindent \textbf{Step 3}.
We may assume (up to a subsequence) that  $\vgf_k \rightharpoonup u$ in $W^{1,p}_{\loc}(\Omega)$, and $\vgf_k \to u$ in $L_{\loc}^{p}(\Omega)$. Let $\omega \Subset \Omega$ be a smooth domain, and denote
$$\mathcal{Q}^{\omega}_{A}(u):=\int_{\omega}L_{A}(u,v)\dx=\int_{\omega}\mathcal{L}_A(x,u, \nabla u)\dx \qquad u\in W^{1,p}(\omega).$$
\noindent {\bf Claim:} $\mathcal{Q}^{\omega}_{A}(u)$ is weakly lower semicontinuous in $W^{1,p}(\omega)$.
Indeed, the sum of the first two terms of the Lagrangian $\mathcal{L}_A(x,z, q)$ is equal to  $|q|^{p}_{A}+(p-1)\frac{|\nabla v|^{p}_{A}}{v^{p}}|z|^{p}$ which is convex in $q$. Therefore, the corresponding functional is weakly lower semicontinuous in $W^{1,p}(\omega)$. So, it remains to prove that the functional
\begin{equation} \
\mathcal{J}^{\omega}_{A}(u):=\int_{\omega}\frac{u^{p-1}}{v^{p-1}}\nabla u\cdot  A(x)\nabla v|\nabla v|^{p-2}_{A}\dx \qquad u\in W^{1,p}(\omega)
\end{equation}
is weakly continuous on any sequence $\left\{\vgf_k\right\}$ satisfying $\vgf_k \!\rightharpoonup \!u$ in $W^{1,p}(\omega)$. Indeed
\begin{multline} \label{55}
\mathcal{J}^{\omega}_{A}(\vgf_k)-\mathcal{J}^{\omega}_{A}(u)=\int_{\omega}\frac{|\nabla v|^{p-2}_{A}}{v^{p-1}}\nabla v\cdot  A(x) \nabla \vgf_k\left(\vgf^{p-1}_{k}-u^{p-1}\right)\dx+ \\ \int_{\omega}\nabla (\vgf_k-u) \cdot A(x)\nabla v|\nabla v|^{p-2}_{A}\frac{u^{p-1}}{v^{p-1}}\dx.
\end{multline}
For a renamed subsequence, there exists a $U \in L^{p}(\omega)$, such that $0 \leq \vgf_k \leq U$ and $\vgf_k \rightarrow u$ a.e. in $\omega$. Therefore $\vgf^{p-1}_{k} \leq U^{p-1} \in L^{P'}(\omega)$. Since $\vgf_k, u, U$ are nonnegative we have that
\begin{equation} \label{56}
\left|\vgf^{p-1}_{k}-u^{p-1}\right|^{p'} \leq \left(U^{p-1}+u^{p-1}\right)^{p'} \leq C\left|U^{p}+u^{p}\right| \in L^{1}(\omega).
\end{equation}
Hence by H\"{o}lder's inequality and Lebesgue's dominated convergence theorem,
\begin{equation}
\left|\int_{\omega}\frac{1}{v^{p-1}}|\nabla v|^{p-2}_{A}\nabla v A(x) \nabla \vgf_k\left(\vgf^{p-1}_{k}-u^{p-1}\right)\dx\right| \leq  C\left\|\nabla \vgf_k\right\|_{p} \left\|\vgf^{p-1}_{k}-u^{p-1}\right\|_{p'} \rightarrow 0.
\end{equation}

Consider the functional $$\Phi(w):=\int_{\omega}\nabla w \cdot A(x)\nabla v|\nabla v|^{p-2}_{A}\frac{u^{p-1}}{v^{p-1}}\dx \qquad w\in W^{1,p}(\omega).$$
Note $A(x)\nabla v|\nabla v|^{p-2}_{A}v^{1-p}u^{p-1} \in L^{p'}(\omega)$. Therefore, H\"{o}lder's inequality implies that $\Phi$ is a continuous functional on $W^{1,p}(\omega)$. Hence, by the definition of weak convergence, the second term of the right hand side of (\ref{55}) converges to zero.

Therefore, $\mathcal{Q}^{\omega}_{A}(u)$ is weakly lower semicontinuous, and we conclude that $0 \leq \mathcal{Q}^{\omega}_{A}(u) \leq \liminf_{k\to\infty} \mathcal{Q}^{\omega}_{A}(\vgf_k)=0$. Moreover, $\int_{B}u^{p}\dx=1$ and $u$ is a ground state. Since $\mathcal{Q}^{\omega}_{A}(u)=0$ for every subdomain $\omega \Subset \Omega$ containing $B$, it follows that $L_{A}(u,v)=0$, and by Proposition~\ref{Picone} $u=cv$, where  $c^{-1}=\big(\int_{B}v^{p}\dx\big)^{1/p}$. Hence, any null sequence converges to the same positive supersolution $v$, and this implies the uniqueness claim.
\end{proof}

We prove now Theorem~\ref{main thm}.
\begin{proof}[Proof of Theorem~\ref{main thm}] Lemma~\ref{crit} implies part~(1) and the necessity parts of (2) and (3).
On the other hand if $\mathcal{Q}_{A,V}$  admits a null sequence, then by Lemma~\ref{crit},  (\ref{derivative}) admits a unique positive continuous supersolution. The latter implies
that the inequality $\mathcal{Q}_{A,V}\geq 0$ cannot be improved, Consequently, $\mathcal{Q}_{A,V}$ is critical in $\Gw$, and hence null sequences exist.

To complete the proof of part (2), we need to show that if $\mathcal{Q}_{A,V}$  is critical in $\Gw$, then it admits also a null sequence that converges locally uniformly in $\Omega$.
Let $\left\{\Omega_{N}\right\}^{\infty}_{N=1}$ be an exhaustion of $\Omega$ such that $x_{0} \in \Omega_{1}$. Pick a nonzero nonnegative function $W \in C^{\infty}_{0}(\Omega_{1})$. For $t \geq 0$ and $N \geq 1$, consider the functional $\mathcal{Q}_{A,V-tW}$ on $C^{\infty}_{0}(\Omega_{N})$. By Proposition~\ref{pro4}, there exists a unique $t_{N}>0$ such that the functional $\mathcal{Q}_{A,V-t_{N}W}$ is critical in $\Omega_{N}$ (note the proof of Proposition~\ref{pro4} relies only on lemmas \ref{subcrit} and \ref{crit} but not on our theorem). Denote by $\phi_{N}$ the corresponding ground state satisfying $\phi_{N}(x_{0})=1$.  Clearly, $\{t_N\}$ is a positive nonincreasing sequence that converges to $t_\infty \geq 0$. By Harnack's convergence principle, up to a subsequence, $\{\phi_N\}$ converges locally uniformly to a positive solution $v$ of the equation $Q_{A,V-t_\infty W}(u)=0$ in $\Gw$.

On the other hand, since $\mathcal{Q}_{A,V}$ is critical in $\Gw$, it follows that for any $t>0$ the functional $\mathcal{Q}_{A,V-tW} \not\geq 0$ on $\core$, hence $t_{\infty}= 0$. By the uniqueness of the positive (super)solutions of the equation $Q_{A,V}(u)=0$ in $\Gw$,  it follows that the whole sequence converges locally uniformly to $v$. Since $\gl_1(Q_{A,V-t_{N}W},\Gw_N) =0$, it follows that  $\phi_{N}$ is a principal eigenfunction of $\mathcal{Q}_{A,V-t_{N}W}$ and therefore, $\phi_{N}\in W^{1,p}_0(\Gw_N)\subset W^{1,p}_0(\Gw)$. Consequently,
$$\mathcal{Q}_{A,V}(\phi_{N})=t_{N}\int_{\Omega}W|\phi_{N}|^p\dx \to 0, \qquad \mathrm{and} \  \int_{B}|\phi_{N}|^{p}\dx \asymp 1.$$
By Remark~\ref{ground}, $\left\{\phi_{N}\right\}$ is a null sequence of $\mathcal{Q}_{A,V}$ and $v$ is the ground state of $\mathcal{Q}_{A,V}$.

In light of Lemma~\ref{crit}, part (4) is proved in Lemma~\ref{subcrit}. So, it remains to prove part~(5). Consider the functional
$$\tilde{\mathcal{Q}}(\vgf):=\mathcal{Q}_{A,V}(\vgf)+\int_{B}|\vgf|^{p}\dx \qquad   \vgf\in \core,$$
where $B\Subset \Omega,$ is a fixed open set. Clearly, $\tilde{\mathcal{Q}}$ is subcritical in $\Gw$. By part~(4), we have that for any $B\Subset \Omega$, there exists a positive continuous function $W$ in $\Omega$ such that
\begin{equation}\label{eq_70}
\int_{\Omega}W(x)|\vgf(x)|^{p}\dx \leq \mathcal{Q}_{A,V}(\vgf)+\int_{B}|\vgf|^{p}\dx \qquad\forall \vgf \in C^{\infty}_{0}(\Omega).
\end{equation}
Suppose now that for each open ball $B\Subset \Omega$ there is a nonnegative sequence $\left\{\vgf_k\right\} \subset C^{\infty}_{0}(\Omega)$ such that $\int_{B}|\vgf_k|^{p}\dx=1, \mathcal{Q}_{A,V}(\vgf_k) \to 0$ and $\int_{\Omega}\vgf_k\psi\dx \to 0$. Since $\left\{\vgf_k\right\}$ is a null sequence, by Lemma~\ref{crit}, it converges to a ground state $\phi>0$ in $L^{p}_{\loc}(\Omega)$. Then $\int_{\Omega}\vgf_k\psi\dx \to \int_{\Omega}\phi\,\psi\dx \neq 0$, and we arrive at a contradiction.
Therefore, there exists an open $B\Subset \Omega$ such that
\begin{equation}\label{eq_71}
\int_{B}|\vgf|^{p}\dx \leq C \left(\mathcal{Q}_{A,V}(\vgf)+\left|\int_{\Omega}\vgf\,\psi\dx\right|^{p}\right) \qquad \forall \vgf \in C^{\infty}_{0}(\Omega).
\end{equation}
Hence, \eqref{eq_70} and \eqref{eq_71} imply \eqref{poincare}.
\end{proof}
%%%%%%%%%%%%%%%%%
As a consequence we generalize Proposition 2 in \cite{AA} concerning the principle eigenvalue.
%%%%%%%%%%%%%%%%%%%%%%%%%%%%%%%
\begin{lemma} \label{mes lemma}
Let $\Gw$ be a bounded domain, $A\in L^\infty(\Gw,\R^{n^2})$ symmetric and uniformly positive definite, $V\in L^\infty(\Gw)$, and $\mathcal{Q}_{A,V}\geq 0$ in $\Gw$. If $p<2$ assume further that $A\in C^\ga(\Gw)$.
\begin{enumerate}
  \item If (\ref{ev.pr}) admits a real eigenvalue $\lambda$ with an eigenfunction $\psi\geq 0$, then $\gl=\gl_1$.
  \item $\gl_1$ is a simple eigenvalue.
  \item The principal eigenfunction associated with $\gl_1$ is the unique positive supersolution of the equation $Q_{A,V-\gl_{1}}(u)=0$ in $\Gw$.
    \item If $\lambda\neq\lambda_{1}$ is a real eigenvalue with an eigenfunction $\psi$, then
\begin{equation} \label{mes}
\mathrm{meas}(\Gw^{-}) \geq \left(\left\|\lambda-V\right\|_{\infty}\bar{C}^{p}\right)^{\sigma},
\end{equation}
where $\Gw^{-}=\{x \in \Gw : \psi(x)<0\}$, $\bar{C}$ is a constant independent of $\psi$ and $\lambda$, and $\sigma :=-\frac{n}{p}$ if $p<n$, $\gs:=-2$ if $p \geq n$.

  \item $\gl_1$ is an isolated eigenvalue in $\R$.
\end{enumerate}
 \end{lemma}
%%%%%%%%%%%%%%%%%%%
\begin{Rem} \label{e.v. neq V} {\em
If $\gl$ is an eigenvalue of \eqref{ev.pr} and $\Gw\Subset \R^n$, then $\|\lambda-V\|_{\infty}> 0$. Otherwise, an associated eigenfunction $u$ is  $p\,$-harmonic in $\Gw$ and vanishes on $\partial\Omega$. Hence, $u=0$ in $\Omega$.
 }
\end{Rem}
%%%%%%%%%%%%%%%%%%%%%%%%%%%%%
\begin{proof}[Proof of Lemma~\ref{mes lemma}]
1--3. Any eigenfunction with eigenvalue $\gl_1$ is a minimizer of \eqref{p.ev}. Lemma~\ref{pos.e.f.} implies that \eqref{p.ev}
admits a minimizer $\phi$, and that each minimizer is a principal eigenfunction of the operator $Q_{A,V}$ with eigenvalue $\gl_1$ (up to a sign change).
 Moreover, if $\psi$ is a principal eigenfunction with eigenvalue $\gl$, then
$\psi$ is a positive supersolution of the equation  $Q_{A,V-\gl_{1}}(u)=0$ in $\Gw$.
Since any minimizing sequence is a null sequence (up to a sign change), part (3) of Theorem~\ref{main thm} implies that $\gl=\gl_1$, that $\gl_1$ is simple, and that
$\phi$ is the unique positive supersolution of the equation  $Q_{A,V-\gl_{1}}(u)=0$ in $\Gw$.

4. Recall that $|\psi_{-}| \in W^{1,p}_{0}(\Omega)$, and by Lemma~\ref{subsol}, $|\psi_{-}|$ is a subsolution of the equation $Q_{A,V-\gl}(u)=0$ in $\Gw$. Hence,
$$\left\|\nabla \psi_{-}\right\|^{p}_{p,A} \leq \int_{\Omega}(\lambda-V(x))|\psi_{-}|^{p}\dx\leq \left\|\lambda-V\right\|_{\infty} \int_{\Omega}|\psi_{-}|^{p}\dx.$$
H\"{o}lder inequality implies that
\begin{equation}
\left\|\nabla \psi_{-}\right\|^{p}_{p,A} \leq \left\|\lambda-V\right\|_{\infty}\left\|\psi_{-}\right\|^p_{p^{*}}\left(\mathrm{meas}(\Omega_{-})\right)^{1-\frac{p}{p^{*}}},
\end{equation}
where $p^{*}\!=\!\frac{pn}{n\!-\!p}$ if $p<n$, and $p^{*}\!=\!2p$ if $p\! \geq\! n$.
By Gagliardo-Nirenberg-Sobolev inequality,
$$\left\|\psi_{-}\right\|_{p^{*}} \leq C\left\|\nabla \psi_{-}\right\|_{p} \leq \bar{C}\left\|\nabla \psi_{-}\right\|_{p,A}.$$
Therefore,
$$\left\|\nabla \psi_{-}\right\|^{p}_{p,A} \leq \left\|\lambda-V\right\|_{\infty}\bar{C}^{p}\left\|\nabla \psi_{-}\right\|^{p}_{p,A}\left(\mathrm{meas}(\Omega_{-})\right)^{1-\frac{p}{p^{*}}},$$
and we obtain the desired result.

5. Assume that there exists a sequence of eigenvalues $\{\lambda_{k}\}\subset \R$ such that $\lambda_{k}>\lambda_{1}$, $\lim_{k \to \infty}\lambda_{k}=\lambda_{1}$. Let $\{\phi_{k}\} \in W_{0}^{1,p}(\Omega)$ be the corresponding normalized eigenfunctions. By compactness and the simplicity of $\lambda_{1}$,  $\lim_{k \to \infty}\phi_{k}=\phi$ in $W_{0}^{1,p}(\Omega)$, where $\phi>0$ is the principal eigenfunction with eigenvalue $\lambda_{1}$. By part~(4), and Remark~\ref{e.v. neq V},  $\mathrm{meas} \{\phi_{k}<0\}\geq c >0$, where $c$ is independent of $k$, a contradiction to the positivity of $\phi$.
\end{proof}
\begin{corollary}\label{cor_simple1}
Assume that $\Gw$ is a bounded domain,  $A$ is a bounded measurable symmetric matrix which is uniformly positive definite in $\Omega$, and $V\in L^\infty(\Gw)$. If $p<2$ assume further that $A\in C^\ga(\Gw)$. Then all the assertions of Theorem~\ref{Weak, Strong, Pos.} are equivalent.
\end{corollary}
\begin{proof}
$(iv')\Rightarrow (iii)$. By the AAP theorem $\gl_1\geq 0$. Lemma~\ref{mes lemma} implies that the principal eigenfunction is the unique positive supersolution of the equation $Q_{A,V-\gl_{1}}(u)=0$ in $\Gw$. Hence, $\gl_1>0$ and $(iv')\Rightarrow (iii)$  (cf. the proof of this part in \cite{GS}, where the WCP was used under the assumption $\gl_1\leq 0$!). Thus, all the assertions of Theorem~\ref{Weak, Strong, Pos.} are equivalent.
\end{proof}
%%%%%%%%%%%%%%%
As in \cite{PT3}, criticality is characterized in terms of the relative capacity of compact sets.
\begin{definition}{\em
Suppose that the functional $\mathcal{Q}_{A,V}$ is nonnegative on $C^{\infty}_{0}(\Omega)$. Let $K \Subset \Omega$ be a compact set. The {\em $\mathcal{Q}_{A,V}$-capacity} of $K$ in $\Omega$ is defined by
$$\mathrm{Cap}_{\mathcal{Q}_{A,V}}(K,\Omega):= \inf\left\{\mathcal{Q}_{A,V}(\vgf) \mid \vgf \in C^{\infty}_{0}(\Omega), \ \vgf \geq 1 \ \mathrm{on} \ K\right\}.$$
}
\end{definition}
The following theorem generalizes Theorem~3.4 in \cite{PT2} and Theorem~4.5 in \cite{PT3}. We omit the proof since it differs only slightly from the proofs in \cite{PT2,PT3}.
\begin{theorem} \label{W thm}
Suppose that $\mathcal{Q}_{A,V}\geq 0$  in $\Omega$, where  $A\in C^\ga(\Gw,\R^{n^2})$ satisfy conditions (A) and (E), and $V\in L^\infty_\loc(\Gw)$.
Then following statement are equivalent.
\begin{itemize}
	\item [(i)] $\mathcal{Q}_{A,V}$ is subcritical in $\Omega$.
	\item[(ii)] There exists a continuous function $W>0$ in $\Omega$ such that
\begin{equation}\label{57}
\mathcal{Q}_{A,V}(\vgf) \geq \int_{\Omega}W(x)(|\nabla \vgf|^{p}_{A}+|\vgf|^{p})\dx \qquad\forall \vgf\in C_{0}^{\infty}(\Omega).	
\end{equation}
  \item[(iii)]  There exists an open set $B \Subset \Omega$ and $c_{B}>0$ such that
\begin{equation}\label{58}
\mathcal{Q}_{A,V}(\vgf) \geq C_{B}\left|\int_{B}|\vgf|\dx\right|^{p}\qquad\forall \vgf\in C_{0}^{\infty}(\Omega).	
\end{equation}
  \item[(iv)] The $\mathcal{Q}_{A,V}$-capacity of any closed ball in $\Omega$ is positive.
\end{itemize}

Suppose further that $p<n$ and let $p^{*}:=\frac{p\,n}{n-p}$ be its critical Sobolev exponent. Then $\mathcal{Q}_{A,V}$ is subcritical in $\Omega$ if and only if the following  Hardy-Sobolev-Maz'ya-type inequality holds true: there exists a continuous function $\tilde{W}>0$ in $\Omega$ such that
\begin{equation}\label{59}
\mathcal{Q}_{A,V}(\vgf) \geq \left(\int_{\Omega}\tilde{W}(x)|\vgf|^{p^{*}}\dx\right)^{\frac{p}{p^{*}}}\quad\forall \vgf\in C_{0}^{\infty}(\Omega).	
\end{equation}
\end{theorem}

%%%%%%%%%%%%%%%%%%%%%%%%%%%%%%%%%%%%%%%%%%%%%%%%%%%%%%%%%%%%%
\mysection{Criticality Theory}\label{sec_crit_theory}
In this short section we list further positivity properties of $\mathcal{Q}_{A,V}$, generalizing the results of \cite[Section~4]{PT1}.
We omit the proofs since they differ only slightly from the proofs in \cite{PT1}.
\begin{proposition} \label{pro1}
Let $V_{1},V_{2} \in L^{\infty}_{\loc}(\Omega)$ and suppose that $V_{2} \gneqq V_{1}$.

1. If $\mathcal{Q}_{A,V_{1}} \geq 0$ on $C_{0}^{\infty}(\Omega)$, then $\mathcal{Q}_{A,V_{2}}$ is subcritical in $\Gw$.

2. If $\mathcal{Q}_{A,V_{2}}$ is critical in $\Omega$, then $\mathcal{Q}_{A,V_{1}}$ is supercritical in $\Omega$.
\end{proposition}
%%%%%%%%%%%%%%%%%%%%%%%
\begin{proposition} \label{pro2}
Let $\Omega_{1} \subset \Omega_{2}$ be domains in $\Real ^{n}$ such that $\Omega_{2} \setminus \bar{\Omega}_{1} \neq \emptyset$.

1. If $\mathcal{Q}_{A,V} \geq 0$ on $C^{\infty}_{0}(\Omega_{2})$, then $\mathcal{Q}_{A,V}$ is subcritical in $\Omega_{1}$.

2. If $\mathcal{Q}_{A,V}$ is critical in $\Omega_{1}$, then $\mathcal{Q}_{A,V}$ is supercritical in $\Omega_{2}$.
\end{proposition}
%%%%%%%%%%%%%%%%%%%%%%%%%%%%%%%%%%
\begin{proposition} \label{convex}
Let $V_{0},V_{1} \in L^{\infty}_{\loc}(\Omega)$. For $t \in \Real$ we denote
\begin{equation}
\mathcal{Q}_{A,t}(\vgf):=(1-t)\mathcal{Q}_{A,V_{0}}(\vgf)+t\mathcal{Q}_{A,V_{1}}(\vgf) \qquad \vgf\in \core.
\end{equation}
Suppose that $\mathcal{Q}_{A,V_{i}} \geq 0$ on $C^{\infty}_{0}(\Omega)$ for $i=0,1$.
Then the functional $\mathcal{Q}_{A,t} \geq 0$ on $C^{\infty}_{0}(\Omega)$ for all $t \in \left[0,1\right]$.
Moreover, if $V_{0} \neq V_{1}$, then $\mathcal{Q}_{A,t}$ is subcritical in $\Omega$ for all $t \in (0,1)$.
\end{proposition}
%%%%%%%%%%%%%%%%%%%%%%
\begin{proposition} \label{pro4}
Let $\mathcal{Q}_{A,V}$ be a subcritical functional in $\Omega$. Consider $V_{0} \in L^{\infty}(\Omega)$ such that $V_{0} \ngeq 0$ and $\mathrm{supp}V_{0} \Subset \Omega$. Then there exist $\tau_{+}>0$ and $-\infty \leq \tau_{-}<0$ such that $\mathcal{Q}_{A,V+tV_{0}}$ is subcritical in $\Omega$ if and only if $\ t \in (\tau_{-}, \tau_{+})$, and $\mathcal{Q}_{A,V+\tau_{+}V_{0}}$ is critical in $\Omega$.
\end{proposition}
\begin{proposition} \label{pro5}
Let $\mathcal{Q}_{A,V}$ be a critical functional in $\Omega$, and let $\phi$ be the corresponding ground state, Consider $V_{0} \in L^{\infty}(\Omega)$ such that $\mathrm{supp}V_{0} \Subset \Omega$. Then there exists $0< \tau_{+} \leq \infty$ such that $\mathcal{Q}_{A,V+tV_{0}}$ is subcritical in $\Gw$ for $t \in (0,\tau_{+})$ if and only if
$\int_{\Omega}V_{0}|\phi|^{p}\dx>0$.
\end{proposition}
%%%%%%%%%%%%%%%%
\mysection{Liouville Theorem}\label{sec_Liouville}

In the present section we generalize Liouville-comparison theorems proved in \cite{P,PTT} (see also references therein). We note that the results in \cite{PTT} are the counterpart of the results in \cite{P} proved for the linear case ($p=2$) for operators of the form $Pu:=-\nabla \cdot (A(x)\nabla u)+ V(x)u $ with a measurable symmetric matrix valued function $A$ which is locally bounded and locally uniformly positive definite, while \cite{PTT} deals with $p$-Laplace type equation. We have:

\begin{theorem} \label{liouville thm}
For $j=0,1$,  consider the functional
\begin{equation} \label{90}
\mathcal{Q}_{A_{j},V_{j}}(\vgf):=\int_{\Omega}\left(|\nabla \vgf|^{p}_{A_{j}}+V_{j}(x)|\vgf|^{p}\right)\dx  \qquad \vgf \in C^{\infty}_{0}(\Omega),
\end{equation}
where the matrix $A_j$ and the potential $V_{j}$ satisfy assumptions (A), (E), and (V). If $p<2$ assume further that $A_j\in C^\ga(\Gw)$.
Suppose that the following assumptions hold true:
\begin{itemize}
\item[(i)] The functional $\mathcal{Q}_{A_{1},V_{1}}$ admits a ground state $\phi$ in $\Omega$.
\item[(ii)] $\mathcal{Q}_{A_{0},V_{0}} \!\geq \!0$ in $\Omega$, and the equation $Q_{A_{0},V_{0}}(u)\!=\!0$ in $\Omega$ admits a subsolution $\psi \!\in\! W^{1,p}_{\loc}(\Omega)$ satisfying $\psi_{+}\!:=\!\max\!\left\{0,\psi\right\}\! \neq\! 0$.
\item[(iii)] There exists $M > 0$ such that
\begin{equation} \label{91}
\psi_{+}^{2}A_{0} \leq M^2\phi^2A_{1}\qquad \mbox{a.~e. in } \Omega.
\end{equation}
That is,  for almost all $x\in \Omega$ the matrix $\big(M\phi(x)\big)^2A_{1}(x)-\big(\psi_{+}(x)\big)^{2}A_{0}(x)$ is non\linebreak negative-definite on $\mathbb{R}^n$.

\item[(iv)] There exists $N > 0$ such that
\begin{equation} \label{92}
|\nabla \psi|^{p-2}_{A_{0}} \leq N^{p-2}|\nabla \phi|^{p-2}_{A_{1}} \qquad \mbox{a.~e. in } \Omega \cap \left\{\psi > 0\right\}.
\end{equation}
\end{itemize}
Then the functional $\mathcal{Q}_{A_{0},V_{0}}$ is critical in $\Omega$, and $\psi$ is the unique positive supersolution of the equation $Q_{A_{0},V_{0}}(u)=0$ in $\Omega$.
\end{theorem}
%%%%%%%%%%%%%%%%%%%%%%%%%%%%%%%
\begin{proof}[Proof of Theorem~\ref{liouville thm}]
By Lemma~\ref{subsol}, we may assume that $\psi \geq 0$.

\noindent Let $\left\{\vgf_k\right\}\!\subset\!\core$ be a null sequence for $\mathcal{Q}_{A_{1},V_{1}}$. So, there exists an open set $B\! \Subset\! \Omega$ such that $\int_{B}{|\vgf_k|^{p}\dx}=1$, and
$\lim_{k \to \infty}\mathcal{Q}_{A_{1},V_{1}}(\vgf_k)=0$. Without loss of generality, we may assume that $B \subset \mathrm{supp}\, \psi$. Let $w_{k}:=\frac{\vgf_k}{\phi}\geq 0$.  From \eqref{3} it follows that there exists $C_1>0$ such that
\begin{equation}\label{94}
\int_{\Omega}\phi^{2}|\nabla w_{k}|^{2}_{A_{1}}\Big(w_{k}|\nabla \phi|_{A_{1}}+\phi|\nabla w_{k}|_{A_{1}}\Big)^{p-2}\dx\leq C_1\mathcal{Q}_{A_{1},V_{1}}(\vgf_k) \to 0 \quad \mbox{as } k\to\infty.
\end{equation}
Fix $\alpha \in \Real_{+}$, and consider the function $f:\Real^{2}_{+} \to \Real_{+}$ defined by $$f(s,t):=t^{2}\left(\alpha s^{\frac{1}{p-2}}+t\right)^{p-2}.$$
Since
$$f_{s}(s,t)=\alpha t^{2}s^{\frac{3-p}{p-2}}\left(\alpha s^{\frac{1}{p-2}}+ t\right)^{p-3} \geq 0, \quad f_{t}(s,t)=t\left(\alpha s^{\frac{1}{p-2}}+ t\right)^{p-3}\left(2\alpha s^{\frac{1}{p-2}}+ tp\right) \geq 0,$$
the function $f(s,t)$ is a nondecreasing monotone function in each variable separately. Set $t_{0}:=\psi|\nabla w_{k}|_{A_{0}}$ and $t_{1}:=M\phi|\nabla w_{k}|_{A_{1}}$, then assumptions \eqref{91} implies that $t_{0} \leq t_{1}$. Hence,
\begin{equation*}
I_k\!:=\!\!\int_{\Omega}\!\!\psi^{2}|\nabla w_{k}|^{2}_{A_{0}}\!\Big(w_{k}|\nabla \psi|_{A_{0}}\!+\psi|\nabla w_{k}|_{A_{0}}\!\Big)^{p-2}\!\!\!\!\dx \!\leq\!
M^2\!\!\!\int_{\Omega}\!\!\!\phi^{2}|\nabla w_{k}|^{2}_{A_{1}}\!\Big(w_{k}|\nabla \psi|_{A_{1}}\!+M\phi|\nabla w_{k}|_{A_{1}}\!\Big)^{p-2}\!\!\!\!\dx.
\end{equation*}
Let now $s_{0}:=|\nabla \psi|^{p-2}_{A_{0}}$ and $s_{1}:=N^{p-2}|\nabla \phi|^{p-2}_{A_{1}}$, then \eqref{92} implies that $s_{0} \leq s_{1}$. Therefore,
\begin{multline*}
I_k\leq  M^{2}\int_{\Omega}\phi^{2}|\nabla w_{k}|^{2}_{A_{1}}\Big(Nw_{k}|\nabla \phi|_{A_{1}}+M\phi|\nabla w_{k}|_{A_{1}}\Big)^{p-2}\dx\leq\\
C_2\!\!\int_{\Omega}\!\!\phi^{2}|\nabla w_{k}|^{2}_{A_{1}}\Big(w_{k}|\nabla \phi|_{A_{1}}+\phi|\nabla w_{k}|_{A_{1}}\Big)^{p-2}\dx\leq
 C_{2}C_1\mathcal{Q}_{A_{1},V_{1}}(\vgf_k) \to 0.
\end{multline*}
By  \eqref{5} we have $\mathcal{Q}_{A_{0},V_{0}}(\psi w_k) \leq C_3 I_k$. So, $\mathcal{Q}_{A_{0},V_{0}}(\psi w_{k}) \to 0$. Since $w_{k} \to 1$ in $L^{p}_{\loc}(\Omega)$, it follows that $\psi w_{k} \to \psi$ in $L^{p}_{\loc}(\Omega)$. Consequently, $\int_{B}\phi^{p}w_{k}^{p}\dx=1$ implies that $\int_{B}\psi^{p}w_{k}^{p}\dx \asymp 1$. In light of Remark~\ref{ground}, we conclude that $\psi$ is a ground state of $\mathcal{Q}_{A_{0},V_{0}}$.
\end{proof}

\begin{example}[cf. \cite{PTT} examples~3.1--3.3] {\em
Theorem~\ref{liouville thm} implies that we may replace  the $p$-Laplacian appearing as the principal part of  $Q_0$ in examples~3.1--3.3 of \cite{PTT} with the $(p,A_0)$-Laplacian, where $A_0$ is a {\em bounded} measurable symmetric matrix which is locally uniformly positive definite in $\R^n$.
}
\end{example}
\begin{example}\label{ex2} {\em
Let $1 \!\leq \!n \!\leq\! p$, $p\!>\!1$, $\Omega\!=\!\Real^{n}$, and consider the functional
$$\mathcal{Q}_{A_1,V_{1}}(\vgf):=\mathcal{Q}_{I,0}(\vgf)=\int_{\Real^{n}}|\nabla \vgf|^{p}\dx \qquad \vgf\in C_0^\infty(\R^n).$$
By Example~\ref{ex_subcr} the functional $\mathcal{Q}_{A_1,V_{1}}$ admits a ground state $\phi\!=\!\mathrm{const.}$ in $\Real^{n}$.

Consider the functional $$\mathcal{Q}_{A_0,V_0}(\vgf):=\mathcal{Q}_{A_0,0}(\vgf)=\int_{\Omega}|\nabla \vgf|_{A_0}^{p}\dx \qquad \vgf\in C_0^\infty(\R^n),$$
where $A_0$ is  a {\em bounded} measurable symmetric matrix which is locally uniformly positive definite in $\R^n$.  Then $\psi=\mathrm{const.}>0$ is a positive (sub)solution of the equation $Q_{A_0,V_{0}}(u)=0$ in $\Real^{n}$. By Theorem~\ref{liouville thm}, $\psi$ is the unique positive supersolution and unique bounded solution of the equation $Q_{A_0,V_{0}}(u)=0$ in $\Real^{n}$ (cf. 6.10 and 6.11 of \cite{HKM}).
}
\end{example}

%%%%%%%%%%%%%%%%%%%%%%%%%%%%%%
\mysection{Minimal Growth}\label{sec_minimal growth}

In this section we generalize some results in \cite[Section~5]{PT3} concerning the existence of positive solutions of minimal growth in a neighborhood of infinity in $\Omega$.

\begin{definition}{\em
Let $K_{0}$ be a compact set in $\Omega$. A positive solution $u$ of the equation $Q_{A,V}(w)=0$ in $\Omega \setminus K_{0}$ is said to be a \emph{positive solution of minimal growth in a neighborhood of infinity in $\Omega$} (or $u \in \mathcal{M}_{\Omega,K_{0}}$ for brevity) if for any compact set $K$ in $\Omega$, with a smooth boundary, such that $K_{0} \Subset \mathrm{int}(K)$, and any positive supersolution $v \in C((\Omega \setminus K) \cup \partial K)$ of the equation $Q_{A,V}(w)=0$ in $\Omega \setminus K$, the inequality $u \leq v$ on $\partial K$ implies that $u \leq v$ in $\Omega \setminus K$.

 If $u \!\in\! \mathcal{M}_{\Omega,\emptyset}\,$, then $u$ is called a \emph{global minimal solution of the equation $Q_{A,V}(w)=0$ in $\Omega$}.
}
\end{definition}
\begin{Rem} \label{rem1}{\em
Suppose that $\mathcal{Q}_{A,V}\geq 0$ in $\Omega$, and let $v>0$ be a positive solution of \eqref{derivative}. Lemma~\ref{mono} implies that $\lambda_{1}(\tilde{\Omega})>0$ for any bounded domain $\tilde{\Omega} \Subset \Omega$ .

Let $\left\{\Omega_{N}\right\}^{\infty}_{N=1}$ be an exhaustion of $\Omega$. Fix $K \Subset \Omega$ with smooth boundary, and let $u\in C^\ga(\partial K)$ be a positive function. Let $u_{N}$ be a solution of the following Dirichlet problem
\begin{equation} \label{71}
\begin{cases}
Q_{A,V}(w)=0 &\text{ in $ \Omega_{N} \setminus K$},\\
w=u &\text{ on $ \partial K$},\\
w=0 &\text{ on $ \partial \Omega_{N}$}.
\end{cases}
\end{equation}
By Theorem~\ref{Weak, Strong, Pos.}, problem~\eqref{71} is solvable. Moreover, by the WCP (Theorem~\ref{compare thm}), $\{u_N\}$ is a monotone nondecreasing sequence, satisfying $u_N\leq cv$ in $ \Omega_{N} \setminus K$, where $c:=\|\frac{u(x)}{v(x)}\|_{\infty,K}$.
Denote: $u^{K}:=\lim_{N \to \infty}u_{N}$ on  $\Omega \setminus K$.

Clearly, $u^K\leq cv$.  Moreover, a comparison argument implies that $u^K$ does not depend on the exhaustion $\left\{\Omega_{N}\right\}^{\infty}_{N=1}$. Moreover, by the Harnack convergence principle, $u^K$ is a positive solution of the equation  $Q_{A,V}(w)=0$ in $ \Omega \setminus K$.

}
\end{Rem}
A further comparison argument shows that
\begin{lemma} \label{lem1}
 Let $\mathcal{Q}_{A,V}$,  $K$, $v$ and $u$ as above. Then $u^K\leq cv$, and $u^{K} \in \mathcal{M}_{\Omega, K}$.
\end{lemma}

\begin{lemma} \label{lem2}
Let $K_{0} \Subset \Omega$, and let $u$ be a positive solution of the equation $Q_{A,V}(w)=0$ in $\Omega \setminus K_{0}$. Then $u \in \mathcal{M}_{\Omega, K_{0}}$ if and only if for any compact set $K \Subset \Omega$ with smooth boundary, such that $K_{0} \Subset \mathrm{int}(K)$, we have $u=u^{K}$.
\end{lemma}
%%%%%%%%%%%%%%%%%%%%%%%%%%%%%%%%%%%%%%%%
\begin{proof}
Fix a smooth compact set $K \Subset \Omega$  such that $K_{0} \Subset \mathrm{int}(K)$, and let $v \in C((\Omega \setminus K)\cup \partial{K})$ be a positive supersolution of the equation $Q_{A,V}(w)=0$ in $\Omega \setminus K$, satisfying the inequality $u \leq v$ on $\partial{K}$. Then by a standard comparison principle, we conclude that $u^{K}=\lim_{N \to \infty}u_{N} \leq v$ in $\Omega \setminus K$. Hence, if $u=u^{K}$ in $\Omega \setminus K$, it follows that $u \leq v$ in $\Omega \setminus K$. Hence, $u \in \mathcal{M}_{\Omega, K_{0}}$.

On the other hand, if  $u \in \mathcal{M}_{\Omega, K_{0}}$, then $u$ is a positive solution of the equation $Q_{A,V}(w)\!=\!0$ in $\Omega \setminus K_{0}$. Therefore, as above, $u^{K} \leq u$ in $\Omega \setminus K$. By definition, $u \leq u^{K}$ in $\Omega \setminus K$. Thus, $u = u^{K}$ in $\Omega \setminus K$.
\end{proof}
\begin{theorem}
Suppose that $1<p<\infty$, and $\mathcal{Q}_{A,V}\geq 0$ in $\Omega$. Then for any $x_{0} \in \Omega$ the equation $Q_{A,V}(w)=0$ admits a positive solution $u \in \mathcal{M}_{\Omega, \left\{x_{0}\right\}}$.
\end{theorem}

\begin{proof}
Consider an exhaustion $\left\{\Omega_N\right\}^{\infty}_{N=1}$ of $\Omega$ such that $x_{0} \in \Omega_{1}$. Let $\left\{f_{N}\right\}$ be a sequence of nonzero nonnegative smooth functions such that $f_{N}\in C^\infty_0\big(B_{2/N}(x_{0}) \setminus \overline{B_{1/N}(x_{0})}\,\big)$.

Denote $A_{N}:=\Omega_{N} \setminus \overline{B_{1/N}(x_{0})}$,  and choose a fixed reference point $x_{1} \in A_{1}$. Let $v_N$ be a positive solution of the Dirichlet problem
\begin{equation}
\begin{cases}
Q_{A,V}(w)=f_{N} &\text{ in $ A_{N}$},\\
w=0 &\text{ on $ \partial A_{N}$}.
\end{cases}
\end{equation}
Set $u_{N}(x):=v_{N}(x)/v_{N}(x_{1})$. By the Harnack convergence principle, $\left\{u_{N}\right\}$ admits a subsequence which converges locally uniformly in $\Omega \setminus \left\{x_{0}\right\}$ to a positive solution $u$ of the equation $Q_{A,V}(w)=0$ in $\Omega \setminus \left\{x_{0}\right\}$.

Let $K \Subset \Omega$ be a compact set with a smooth boundary such that $x_{0} \in \mathrm{int}(K)$, and let $v \in C((\Omega \setminus K)\cup\partial{K})$ be a positive supersolution of the equation $Q_{A,V}(w)=0$ in $\Omega \setminus K$ satisfying inequality $u \leq v$ on $\partial{K}$.

Now, for $\gd>0$ there exists $N_{K}$ such that  $\mathrm{supp}f_{N}\subset B_{2/N}(x_{0}) \Subset K$ for $N\geq N_K$, and $u_{N} \leq (1+\delta)v$ on $\partial(\Omega_{N} \setminus K)$.
Using a comparison argument, and letting  $N \to \infty$, and then $\delta \to 0$, we obtain that $u \leq v$ in $\Omega \setminus K$. Hence, $u \in \mathcal{M}_{\Omega,\left\{x_{0}\right\}}$.
\end{proof}
%%%%%%%%%%%%%%%%%%%%%%%
\begin{theorem}\label{global thm}
Suppose that the matrix $A$ and the potential $V$ satisfy conditions (A), (E) and (V). If $p<2$ assume further that $A\in C^\ga(\Gw)$. Then $\mathcal{Q}_{A,V}$ is subcritical in $\Omega$ if and only if \eqref{derivative} does not admit a global minimal solution in $\Omega$. In particular, $\phi$ is a ground state of \eqref{derivative} if and only if $\phi$ is a global minimal solution of \eqref{derivative}.
\end{theorem}
\begin{proof}
{\bf  Necessity:} Assume that there exists a global minimal solution $u>0$ of the equation $Q_{A,V}(w)=0$ in $\Omega$, and suppose that $\mathcal{Q}_{A,V}$ is subcritical in $\Omega$. Then there exists a positive supersolution $v$ satisfying $Q_{A,V}(v)\gneqq 0$ in $\Omega$.

Clearly, there exists $\varepsilon>0$ such that $\varepsilon u \leq v$ in $\Omega$. Define $$\varepsilon_{0}:=\max\left\{\varepsilon>0 \mid \varepsilon u\leq v \ \mathrm{in} \ \Omega\right\}.$$
Evidently,  $\varepsilon_{0} u \lneqq v$ in $\Omega$. Therefore, there exist $\delta_{1}, \delta_{2}>0$ and $x_{1}\in \Omega$ such that $$(1+\delta_{1})\varepsilon_{0}u(x) \leq v(x) \qquad x \in B_{\delta_{2}}(x_{1}).$$
Hence, by the definition a global minimal solution it follows that
$$(1+\delta_{1})\varepsilon_{0}u(x) \leq v(x) \qquad x \in \Omega \setminus B_{\delta_{2}}(x_{1}).$$
Consequently, $(1+\delta_{1})\varepsilon_{0}u(x) \leq v(x)$ in $\Omega$, which contradicts the definition of $\varepsilon_{0}$. Thus, subcriticality implies the nonexistence of a global minimal solution in $\Omega$.

\noindent {\bf Sufficiency:} Consider an exhaustion $\left\{\Omega_N\right\}^{\infty}_{N=1}$ of $\Omega$ such that $x_{0} \in \Omega_{1}$, and $x_{1} \in \Omega \setminus \Omega_{1}$. Assume that $\mathcal{Q}_{A,V}$ is critical in $\Omega$, and let $\phi$ be its (unique) ground state satisfying $\phi(x_{1})=1$. We need to prove that $\phi$ is a global minimal solution of the equation $Q_{A,V}(w)=0$ in $\Omega$.

Indeed, fix $i \in \Nat$, and let $f_{i} \in C^{\infty}_{0}(B_{1/i}(x_{0}))$ satisfy $0\leq f_i(x)\leq 1$. For $N \geq 1$, let $u_{N,i}$ be a positive solution of the Dirichlet problem
\begin{equation}
\begin{cases}
Q_{A,V}(w)=f_{i} &\text{in $\Omega_{N}$},\\
w = 0 &\text{on $\partial\Omega_{N}$}.
\end{cases}
\end{equation}
By the WCP,  $\left\{u_{N,i}\right\}_{N \geq 1}$ is a nondecreasing sequence. If $\left\{v_{N.i}(x_{1})\right\}$ is bounded, then the sequence converges to $v_{i}$, where $v_{i}$ satisfies $Q_{A,V}(v_{i})=f_{i} \gneqq 0$ in $\Omega$. Due to Theorem~\ref{main thm} (part~(3)), this contradicts our criticality assumption.  Therefore, $v_{N,i}(x_{1}) \to \infty$ as $N\to \infty$.

Denote $u_{N,i}(x):=\frac{v_{N,i}(x)}{v_{N,i}(x_{1})}$. By Harnack converges principle, we may extract a subsequence of $\left\{u_{N.i}\right\}$ that converges as $N \to \infty$ to a positive solution $u_{i}$ of the equation $Q_{A,V}(w)=0$ in $\Omega$. By the uniqueness of the ground state, we have $u_{i}=\phi$.

Let $K \Subset \Omega$ be a smooth compact set, we may assume that $x_{0} \in \mathrm{int}(K)$. Let $v \in C(\Omega \setminus \mathrm{int}(K))$ be a positive supersolution of the equation $Q_{A,V}(w)=0$ in $\Omega \setminus K$ such that the inequality $\phi \leq v$ holds on $\partial K$. Let $i \in \Nat$ be sufficiently large number such that $\mathrm{supp}f_{i} \Subset K$. For any $\gd>0$ there exists $N_\gd$ such that for $N\geq N_\gd$ we have
\[
\begin{cases}
0=Q_{A,V}(u_{N,i})\leq Q_{A,V}(v) &\text{in $\Omega_{N} \setminus K$},\\
Q_{A,V}(v)\geq 0 &\text{in $\Omega_{N} \setminus K$},\\
0\leq u_{N,i}\leq (1+\gd) v &\text{on $\partial(\Omega_{N} \setminus K)$},\\
\end{cases}
\]
which implies that $\phi=u_{i} \leq (1+\gd) v$ in $\Omega \setminus K$. Letting $\gd\to 0$ we obtain $\phi \leq  v$ in $\Omega \setminus K$. Since $K \Subset \Omega$ is an arbitrary smooth compact set, it follows that the ground state $\phi$ is a global minimal solution of the equation $Q_{A,V}(w)=0$ in $\Omega$.
\end{proof}
Suppose that  $u$ is a positive solution of the equation $Q_{A,V}(w)\!=\!0$ in a punctured neighborhood of $x_{0}$, and $1\!<\!p\!\leq\! n$. Then by \cite{JS,JS1}, either $u$ has a removable singularity, or
\begin{equation}
u(x) \asymp  \begin{cases}
|x-x_{0}|^{\alpha(n,p)} &\text{$p<n$},\\
-\mathrm{log}|x-x_{0}|  &\text{$p=n$},\\
\end{cases}
\qquad \mathrm{as} \ \ x \to x_{0},
\end{equation}
where $\alpha(n,p):=(p-n)/(p-1)$. In particular, in the nonremovable case $\lim_{x \to x_{0}}u(x)=\infty$.

Consequently, we have
\begin{theorem} \label{nonremove thm} Suppose that the matrix $A$ and the potential $V$ satisfy conditions (A), (E) and (V). If $p<2$ assume further that $A\in C^\ga(\Gw)$.
Let $x_{0} \in \Omega$, and let $u \in \mathcal{M}_{\Omega,\left\{x_{0}\right\}}$. Suppose that $\mathcal{Q}_{A,V}$ is subcritical in $\Omega$, then $u$ has a nonremovable singularity at $x_{0}$.

Assume that $1<p\leq n$, $x_{0} \in \Omega$, and $u \in \mathcal{M}_{\Omega,\left\{x_{0}\right\}}$. Suppose that $u$ has a nonremovable singularity at $x_{0}$, then $\mathcal{Q}_{A,V}$ is subcritical in $\Omega$.
\end{theorem}
\begin{proof}
Let $u \!\in \! \mathcal{M}_{\Omega,\left\{x_{0}\right\}}$. If $u$ has a removable singularity at $x_0$, then its continuous extension $\bar{u}$ is a global minimal solution in $\Omega$. Hence, by Theorem~\ref{global thm},  $\mathcal{Q}_{A,V}$ is critical.

Assume that $1<p\leq n$. Let $u \in \mathcal{M}_{\Omega,\left\{x_{0}\right\}}$ and suppose that $u$ has a nonremovable singularity at $x_{0}$. If $\mathcal{Q}_{A,V}$ is critical in $\Omega$, then by Theorem~\ref{global thm}, there exists a global minimal solution $v$ in $\Omega$. Let $\vge>0$. Since $\lim_{x \to x_{0}}u(x)=\infty$, a comparison argument implies that $ v \leq \varepsilon u$ in $\Omega$. Hence, $v=0$ which is a contradiction. Therefore, $\mathcal{Q}_{A,V}$ is subcritical  in $\Omega$.
\end{proof}
\begin{definition}
{\em  A function $u\in \mathcal{M}_{\Omega,\left\{x_{0}\right\}}$ having a nonremovable singularity at $x_{0}$ is called a {\em minimal positive Green function} of
$Q_{A,V}$ in $\Gw$. We denote such a function by $G_{A,V}^\Gw (x,x_0)$.
 }
\end{definition}
\begin{problems}
{\em 1. Prove the uniqueness of the positive minimal Green function $G_{A,V}^\Gw (x,x_0)$ and study its {\em asymptotic} behavior as $x\to x_0$.

2. Assume that $p>n$, and consider a nonnegative functional $\mathcal{Q}_{A,V}$. Is it true that $\mathcal{Q}_{A,V}$ is subcritical in $\Omega$ if $Q_{A,V}$ admits a positive minimal Green function $G_{A,V}^\Gw (x,x_0)$?

Note that \cite{FP1,PT3} give affirmative answers to the above problems for the case $A=I$.
 }
\end{problems}
%%%%%%%%%%%%%%%%%%%%%%
%%%%%%%%%%%%%%
\begin{center}{\bf Acknowledgments}
\end{center}

This work is based on the M.~Sc. thesis of N.~R. that was carried out under the supervision of Y.~P..
The authors wish to thank G.~Lieberman and G.~Psaradakis for valuable discussions.
Y.~P. acknowledges the support of the Israel Science
Foundation (grants No. 963/11) founded by the Israel Academy of
Sciences and Humanities.
N.~R. acknowledges the generous support of the Anne R. and Ned Bord Fellowship, and the Andrew Norman Fried Student Aid Fund.

%%%%%%%%%%%%%%%%%%%%%%%%%%%%%%%%%%%%%%%%%%


\begin{thebibliography}{1}

%\newcommand{\cs}{\vspace{2\parskip}}

\bibitem{AA} A.~Anane, Simplicit\'e et isolation de la premi\'ere valeur propre du $p$-Laplacian avec poides,
\emph{C.~R. Acad. Sci. Paris Ser. I Math.} \textbf{305} (1987), 725--728.

\bibitem{AH1} W.~Allegretto, and Y.~X.~Huang, A Picone's identity for the $p$-Laplacian and applications, \emph{Nonlinear Anal.} \textbf{32} (1998), 819--830.

\bibitem{AH2} W.~Allegretto, and Y.~X.~Huang, Principal eigenvalues and Sturm comprison via Picone's identity, \emph{J.~Differential Equations} \textbf{156} (1999), 427--438.

\bibitem{CFKS} H.~L.~Cycon, R.~G.~Froese, W.~Kirsch, and B.~Simon, ``Schr\"{o}dinger Operators with Applications to Quantum Mechanics and Global Geometry'', Texts and Monographs in Physics, Springer Verlag, Berlin, 1987.

\bibitem{DJ} J.~I.~D\'iaz, ``Nonlinear Partial Differential Equations and Free Boundaries'', Pitman Advanced Publishing Program, Vol. \textbf{106}, Boston-London-Malbourne, 1985.

\bibitem{DS} J.~I.~D\'iaz, and E.~Saa, Existence et unicit\'e de solutions positives pour certioanes \'equations elliptiques quasilin\'eaires, \emph{C. R. Acad. Sci. Paris Ser. I Math.} \textbf{305} (1987), 521--524.

\bibitem{FP1} M.~Fraas,  and Y.~Pinchover, Positive Liouville theorems and asymptotic behavior for $p$-Laplacian type elliptic equations with a Fuchsian potential, {\em Confluentes Mathematici} {\bf 3} (2011), 291--323.

\bibitem{FP} M.~Fraas, and Y.~Pinchover, Isolated singularities of positive solutions of
$p$-Laplacian type equations in $\Real^d$, \emph{J.~Differential Equations} \textbf{254} (2013), 1097--1119

\bibitem{GS} J.~Garc\'ia-Meli\'an, and J.~Sabina de Lis, Maximum and comparison principles for operators involving the $p$-Laplacian, \emph{J. Math. Anal. Appl.} \textbf{218} (1998), 49--65.

\bibitem{HKM} J.~Heinonen, T.~Kilpel\"{a}inen, and O.~Martio,
 {\small \it Nonlinear Potential Theory of Degenerate Elliptic Equations},
unabridged republication of the 1993 original, (Dover Publications, Inc., 2006).

\bibitem{LU}  O.~A.~Ladyzhenskaya,  and N.~N.~Ural'tseva, ``Linear and Quasilinear Elliptic Equations", Academic Press, New York-London, 1968.

\bibitem{LL}  E.~H.~Lieb, and M.~Loss,  ``Analysis" (Second edition), Graduate Studies in Mathematics, 14, American Mathematical Society, Providence, RI, 2001.

\bibitem{GML} G.~M.~Lieberman, Boundary regularity for solutions of degenerate elliptic equations, \emph{Nonlinear Anal.} \textbf{12} (1988), 1203--1219.

\bibitem{MP} \'E.~Mitidieri, and S.~I.~Pokhozhaev, Some generalizations of Bernstein's
theorem, \emph{Differ. Uravn}. \textbf{38} (2002), 373--378; translation in \emph{Differ. Equ.}. \textbf{38} (2002), 392--397.

\bibitem{P} Y.~Pinchover, A Liouville-type theorem for Schr\"odinger operators, \emph{Comm. Math. Phys}. \textbf{272} (2007), 75--84.

\bibitem{P07} Y.~Pinchover,  Topics in the theory of positive solutions of second-order elliptic and parabolic partial
differential equations, in ``Spectral Theory and Mathematical Physics: A Festschrift in Honor of Barry Simon's 60th Birthday",
eds. F.~Gesztesy, et al., Proceedings of Symposia in Pure Mathematics {\bf 76}, American Mathematical Society, Providence,
RI, 2007, 329--356.

\bibitem{PTT} Y.~Pinchover, A.~Tertikas, and K.~Tintarev, A Liouville-type theorem for $p$-Laplacian with pootential term, \emph{Ann. Inst. H.~Poincar\'e-Anal. Non Lin\'eaire} \textbf{25} (2008), 357--368.

\bibitem{PT1} Y.~Pinchover, and K.~Tintarev, Ground state alternative for $p$-Laplacian with potential term, \emph{Calc. Var. Partial Differential Equation} \textbf{28} (2007), 179--201.

\bibitem{PT2} Y.~Pinchover, and K.~Tintarev, On positive
solutions of $p$-Laplacian-type equations, in: ``Analysis, Partial
Differential Equations and Applications - The Vladimir Maz'ya Anniversary
Volume''. eds. A.~Cialdea et al., Operator Theory: Advances and Applications,
Vol. \textbf{193}, Birkauser Verlag, Basel (2009), 245--268.

\bibitem{PT3} Y.~Pinchover, and K.~Tintarev, On positive solutions of minimal growth for singular $p$-Laplacian with potential term, \emph{Adv. Nonlinear
Stud.} \textbf{8} (2008) 213--234.

\bibitem{PS} P.~Pucci, and J.~Serrin, ``The Maximum Principle'', Birkh\"{a}user-Verlag, Basel, 2007.

\bibitem{SY} R.~Schoen, and S.-T.~Yau, ``Lectures on Differential Geometry",  Conference Proceedings and Lecture Notes in Geometry and Topology \textbf{I}, International Press, Cambridge, 1994.

\bibitem{JS} J.~Serrin, Local behavior of solutions of quasi-linear elliptic equations, \emph{Acta Math.} \textbf{111} (1964), 247--302.

\bibitem{JS1} J.~Serrin, Isolated singularities of solutions of quasi-linear equations, \emph{Acta Math.} \textbf{113} (1965), 219--240.

\bibitem{Y} N.~Yoshida, ``Oscillation Theory of Partial Differential Equations", World Scientific Publishing Co. Pte. Ltd., Hackensack, NJ, 2008.

\end{thebibliography}
\end{document}